\documentclass[twocolumn]{autart}    

\newtheorem{theorem}{Theorem}
\newtheorem{remark}{Remark}
\newtheorem{lemma}{Lemma}
\newtheorem{proposition}{Proposition}

\newtheorem{assumption}{Assumption}

\usepackage{graphicx}
\usepackage[figuresright]{rotating}
\usepackage{bbm}
\usepackage{graphics}
\usepackage{epsfig}
\usepackage{mathptmx}
\usepackage{times}
\usepackage{amsmath}
\usepackage{amssymb,amsfonts}
\usepackage{float}

\usepackage{cite}
\usepackage{subfigure}
\usepackage{psfrag}
\usepackage{url}
\usepackage{stfloats}
\usepackage[usenames]{color}
\usepackage{algorithm,algorithmicx}

\begin{document}

\begin{frontmatter}

\title{Distributed Algorithms for Computing a Fixed Point of Multi-Agent Nonexpansive Operators\thanksref{footnoteinfo}} 

\thanks[footnoteinfo]{Corresponding author L.~Xie.}

\author[Paestum]{Xiuxian Li}\ead{xiuxianli@ntu.edu.sg},
\author[Paestum]{Lihua Xie}\ead{elhxie@ntu.edu.sg},

\address[Paestum]{School of Electrical and Electronic Engineering, Nanyang Technological University, 50 Nanyang Avenue, Singapore 639798}

\begin{keyword}
Distributed algorithms; multi-agent networks; $Q$-strong-connectivity; nonexpansive operators; fixed point; distributed optimization.
\end{keyword}

\begin{abstract}     
This paper investigates the problem of finding a fixed point for a global nonexpansive operator under time-varying communication graphs in real Hilbert spaces, where the global operator is separable and composed of an aggregate sum of local nonexpansive operators. Each local operator is only privately accessible to each agent, and all agents constitute a network. To seek a fixed point of the global operator, it is indispensable for agents to exchange local information and update their solution cooperatively. To solve the problem, two algorithms are developed, called distributed Krasnosel'ski\u{\i}-Mann (D-KM) and distributed block-coordinate Krasnosel'ski\u{\i}-Mann (D-BKM) iterations, for which the D-BKM iteration is a block-coordinate version of the D-KM iteration in the sense of randomly choosing and computing only one block-coordinate of local operators at each time for each agent. It is shown that the proposed two algorithms can both converge weakly to a fixed point of the global operator. Meanwhile, the designed algorithms are applied to recover the classical distributed gradient descent (DGD) algorithm, devise a new block-coordinate DGD algorithm, handle a distributed shortest distance problem in the Hilbert space for the first time, and solve linear algebraic equations in a novel distributed approach. Finally, the theoretical results are corroborated by a few numerical examples.
\end{abstract}

\end{frontmatter}

\section{Introduction}\label{s1}

Fixed point theory in Hilbert spaces has a wide range of applications, for instance, in optimization, game theory, and nonlinear numerical analysis, and so forth \cite{bauschke2017convex,cegielski2012iterative}. A fixed point for an operator is defined to be a point that maps the point to itself, and an important problem is to develop algorithms to find a fixed point for an operator or to seek a common fixed point for a finite/infinite family of operators.

Up to date, a wide spectrum of works have been established to deal with the fixed point seeking problem for nonexpansive operators under various scenarios \cite{iiduka2016convergence,cominetti2014rate,matsushita2017convergence,bravo2018rates,liang2016convergence,borwein2017convergence,kanzow2017generalized,shehu2018convergence,
dall2019convergence,themelis2019supermann}. It is noteworthy that the conventional Picard iteration is generally not guaranteed to converge for nonexpansive operators, thus engendering a vast amount of research. One of most prominent algorithms for computing a fixed point of a nonexpansive operator is the so-called Krasnosel'ski\u{\i}-Mann (KM) iteration \cite{mann1953mean,krasnosel1955two}, which can converge weakly to a fixed point of the considered nonexpansive operator under mild conditions \cite{reich1979weak}. Along this line, fruitful results have been developed in recent decades. For example, the convergence rate was analyzed for exact and inexact KM iterations in \cite{cominetti2014rate,matsushita2017convergence} and \cite{bravo2018rates,liang2016convergence}, respectively. More recently, a superlinear convergent algorithm, called SuperMann, was proposed and anatomized for computing a fixed point of a nonexpansive operator in \cite{themelis2019supermann}. However, it is worthwhile to point out that all the aforementioned works focus on the centralized case, that is, there exists a global computing unit or coordinator who can access all the information relevant to the considered problem, which is restrictive or even unimplementable for large-scale problems in practice.

Compared with centralized algorithms, distributed algorithms are in possession of a multitude of advantages, including lower computational cost, less storage required, and more robust to failures or malicious attacks, etc. In a distributed scenario, there are a finite number of agents in a network, and each agent is only capable of accessing partial information of the studied problem, hence leading to that every single agent cannot solve the global problem alone. In this case, each agent must propagate its information to some other agents via local communication in order to cope with the global problem in a cooperative manner. Wherein, an agent usually cannot interact with all other agents, instead only with some subset of agents, and the communication channels (e.g., wireless channels) are possibly unreliable and subject to packet dropout, etc., thereby making the communication channels time-varying. Consequently, distributed algorithms under time-invariant/time-varying communication graphs have thus far been vastly investigated in distributed optimization, see \cite{nedic2009distributed,nedic2018multi,liu2017convergence,mateos2017distributed,xie2018distributed,liang2019distributed,li2019distributed,zeng2019generalized,mansoori2019a}, to just name a few, and multi-agent systems/networks, see \cite{olfati2007consensus,ren2010distributed,li2016consensus,meng2017adaptive} for some references. For example, a distributed subgradient algorithm was proposed and analyzed under time-varying interaction graphs in \cite{nedic2009distributed}, while a distributed asynchronous Newton-based algorithm was developed under a static interconnection graph in \cite{mansoori2019a}, showing that the convergence rate can be globally linear and locally quadratic in expectation. Please refer to a recent survey \cite{yang2019survey} for more details on distributed optimization.

Over the past few years, distributed algorithms for finding fixed points of operators have become a focus of researchers. To be specific, the problem of seeking a common fixed point for a finite number of paracontraction operators was addressed in \cite{fullmer2016asynchronous,fullmer2018distributed}, and a distributed algorithm was designed along with convergence analysis. Also, the common fixed point finding problem for a family of strongly quasi-nonexpansive operators was addressed in \cite{liu2017distributed}. The aforementioned common fixed point seeking problem can find numerous applications, such as, in convex feasibility problems \cite{necoara2018randomized,kruger2018set} and solving linear algebraic equations in a distributed approach \cite{mou2015distributed,wang2017further,wang2019distributed,alaviani2018distributed,wang2019scalable}, and so on. It should be noticed that the aforesaid works are only concerned with the Euclidean space. With regard to the Hilbert space, distributed optimization was considered under random communication digraphs in \cite{alaviani2019distributed}, and the common fixed point seeking problem was taken into account for a finite collection of nonexpansive operators in \cite{li2019distributed3}, where two distributed algorithms, namely the distributed inexact KM iteration and distributed inexact block-coordinate KM iteration, were proposed. It can be observed that, in \cite{li2019distributed3}, it is assumed that there exists at least one common fixed point for all local operators, which motivates the further investigation on the case without any common fixed point for all local operators.

Inspired by the above facts, this paper aims at proposing distributed algorithms for finding a fixed point of a global nonexpansive operator under time-varying communication graphs in real Hilbert spaces, where the global operator is separable, consisting of a sum of local nonexpansive operators. In this problem, each local operator is privately known to each agent in a network. To tackle the problem, two distributed algorithms are developed with diminishing stepsizes. In summary, the contributions of this paper are listed below.

\begin{enumerate}
  \item To our best knowledge, this paper is the first to investigate the fixed point seeking problem in real Hilbert spaces in a distributed setup. To cope with the problem, a distributed algorithm, called distributed KM (D-KM) iteration, is developed, and it is shown that, under some appropriate assumptions, the D-KM iteration can converge weakly to a fixed point of the considered global operator.
  \item To alleviate computational complexity, another algorithm, called distributed block-coordinate KM (D-BKM) iteration, is proposed, which is ameliorated based on the D-KM iteration by randomly selecting and computing one block-coordinate of local operators at a point for each agent, instead of the whole coordinate as in the D-KM iteration. It is proved that the D-BKM algorithm is still weakly convergent to a fixed point of the global operator.
  \item The studied problem provides a unified framework for various problems. Specifically, the well-known distributed gradient descent (DGD) algorithm can be recovered from the D-KM iteration, including that a new block-coordinate DGD algorithm is provided by the D-BKM iteration, a distributed shortest distance problem is first resolved in the Hilbert space by appealing to the developed algorithms, and furthermore, the designed D-KM iteration can be leveraged to solve linear algebraic equations in a novel distributed manner.
\end{enumerate}

The rest of this paper is structured as follows. Related basic knowledge is provided in Section II, and the main results are presented in Section III, including two proposed algorithms, i.e., the D-KM and D-BKM iterations. Several applications are given in Section IV, and two numerical examples are introduced in Section V. Finally, the conclusion of this paper is drawn and the future topics are discussed in Section VI.

\section{Preliminaries and problem formulation}\label{s2}

\subsection{Preliminaries}\label{s2.1}

{\em Notations:} Denote by $\mathcal{H}$ a real Hilbert space with inner product $\langle\cdot,\cdot\rangle$ and associated norm $\|\cdot\|$. Define $[N]:=\{1,2,,\ldots,N\}$ to be the index set for any positive integer $N$, and denote by $col(z_1,\ldots,z_k)$ the concatenated column vector or matrix of $z_i,i\in [k]$. For a positive integer $n$, let $\mathbb{R}$, $\mathbb{R}^n$, $\mathbb{R}^{n\times n}$, and $\mathbb{N}$ be the sets of real numbers, $n$-dimensional real vectors, $n\times n$ real matrices, and nonnegative integers, respectively. $P_X(\cdot)$ stands for the projection operator onto a closed and convex set $X\subset\mathcal{H}$, i.e., $P_X(z):=\mathop{\arg\min}_{x\in X}\|z-x\|$ for $z\in\mathcal{H}$. Also, let $I$, $Id$, and $\otimes$ be the identity matrix of appropriate dimension, the identity operator, and the Kronecker product, respectively. Define $d_X(z):=\inf_{x\in X}\|z-x\|$, i.e., the distance from $z\in\mathcal{H}$ to the set $X$. Let $\lfloor c\rfloor$ and $\lceil c\rceil$ be respectively the floor and ceiling functions for a real number $c$. Given an operator $M:\mathcal{H}\to\mathcal{H}$, define $Fix(M):=\{x\in\mathcal{H}:M(x)=x\}$, i.e., the set of fixed points of $M$. Let $\rightharpoonup$ and $\to$ denote weak and strong convergence, respectively, where the weak (resp. strong) convergence to a point $x$ for a sequence $\{x_n\}_{n\in\mathbb{N}}$ is defined as $lim_{n\to\infty}\langle x_n,u\rangle=\langle x,u\rangle$ for all $u\in\mathcal{H}$ (resp. $\lim_{n\to \infty}\|x_n-x\|=0$). Let $B(x;r)$ be the closed ball with center $x$ and radius $r$. Define a $\mathcal{H}$-valued random variable as a measurable map $x:(\Omega,\mathcal{F})\to(\mathcal{H},\mathcal{B})$, along with the probability space $(\Omega,\mathcal{F},\mathbb{P})$ and the expectation $\mathbb{E}$. A measurable (or $\mathcal{F}$-measurable) map is defined by holding $\{\omega\in\Omega: x(w)\in S\}\subset \mathcal{F}$ for every set $S\in\mathcal{B}$. Denote by $\sigma(G)$ the $\sigma$-algebra generated by the family $G$ of random variables. Let $\mathfrak{F}=\{\mathcal{F}_k\}_{k\in\mathbb{N}}$ be a filtration, i.e., each $\mathcal{F}_i$ is a sub-sigma algebra of $\mathcal{F}$ such that $\mathcal{F}_k\subset\mathcal{F}_{k+1}$ for all $k\in\mathbb{N}$. Let $\ell_+(\mathfrak{F})$ be the set of $[0,\infty)$-valued random variable sequence $\{\zeta_k\}_{k\in\mathbb{N}}$ adapted to $\mathfrak{F}$, i.e., $\zeta_k$ is $\mathcal{F}_k$-measurable for all $k\in\mathbb{N}$, and define $\ell_+^1(\mathfrak{F})=\{\{\zeta_k\}_{k\in\mathbb{N}}\in\ell_+(\mathfrak{F}):\sum_{k\in\mathbb{N}}\zeta_k<\infty~a.s.\}$. In this paper, all equalities and inequalities are understood to hold $\mathbb{P}$-almost surely (abbreviated as a.s.) whenever involving random variables, even though ``$\mathbb{P}$-almost surely'' is not expressed explicitly.

In what follows, it is necessary to recall some elementary concepts in operator theory \cite{bauschke2017convex}.

Consider an operator $T:S\to\mathcal{H}$, where $S\subset\mathcal{H}$. The operator $T$ is called {\em nonexpansive} if for all $x,y\in S$,
\begin{align}
\|T(x)-T(y)\|\leq \|x-y\|.         \label{1}
\end{align}

It is well known that the fixed point set $Fix(T)$ is closed and convex if $T$ is nonexpansive \cite{cegielski2015application}.

\subsection{Problem formulation}\label{s2.2}

The objective of this paper is to design distributed algorithms for computing a fixed point for a {\em global operator} $F$ in real Hilbert space $\mathcal{H}$, i.e.,
\begin{align}
\text{find}~x\in\mathcal{H},~~\text{s.t.}~~x\in Fix(F),~~F:=\frac{1}{N}\sum_{i=1}^N F_i,            \label{6}
\end{align}
where $F_i:\mathcal{H}\to\mathcal{H}$ is a nonexpansive operator, called {\em local operator}, for all $i\in[N]$. It can be observed from (\ref{6}) that the global operator $F$ is separable and composed of a sum of local operators, which is especially practical for large-scale problems, and each local operator $F_i$ can be only privately accessible to agent $i$, for example, agent $1$ only knows $F_1$, but unaware of other local operators.

To move forward, let us introduce an assumption, which is needed for the ensuing convergence analysis.

\begin{assumption}\label{a1}
There exists a positive constant $B$ such that
\begin{align}
\|F_i(x)-x\|\leq B,~~~\forall i\in[N],~x\in\mathcal{H}.        \label{7}
\end{align}
\end{assumption}

\begin{remark}\label{rm1}
It should be noted that Assumption \ref{a1} can be satisfied in many interesting problems. For instance, consider the minimization problem for a differentiable convex function $f:\mathbb{R}^n\to\mathbb{R}$ with Lipschitz gradient, that is, $\|\nabla f(x)-\nabla f(y)\|\leq L\|x-y\|$ for all $x,y\in\mathbb{R}^n$ and a constant $L>0$. In this case, the minimization problem is equivalent to the fixed point finding problem for a nonexpansive operator $T:x\mapsto x-\tau\nabla f(x)$, where $\tau\in(0,2/L)$ is a constant \cite{liang2016convergence}. Then, (\ref{7}) in Assumption \ref{a1} amounts to the boundedness of $\nabla f$, which has been widely employed in lots of existing literature in distributed optimization, e.g., \cite{nedic2009distributed,liu2017convergence,mateos2017distributed,xie2018distributed,liang2019distributed}.
\end{remark}

At this position, it is necessary to review graph theory for modeling the communications among all agents \cite{ren2010distributed}. Usually, the $N$ agents construct a network and a digraph $\mathcal{G}=(\mathcal{V},\mathcal{E})$ is exploited to delineate their communications, where $\mathcal{V}=[N]$ and $\mathcal{E}\subset\mathcal{V}\times\mathcal{V}$ are respectively the node (or vertex) and edge sets. An {\em edge} $(i,j)\in\mathcal{E}$ connotes an information flow from agent $i$ to agent $j$, and then agent $i$ (resp. $j$) is called an {\em in-neighbor} or simply {\em neighbor} (resp. out-neighbor) of agent $j$ (resp. $i$). Moreover, a {\em directed path} is a sequence of consecutive edges in the form $(i_1,i_2),(i_2,i_3),\ldots,(i_{l-1},i_{l})$. A graph is called {\em strongly connected} if any node in the graph has at least one directed path to any other node. To be practical, the communication graph is assumed to be time-varying in this paper, i.e., any directed edge can have distinct status at different time slots. As such, let us denote by $\mathcal{G}_k=(\mathcal{V},\mathcal{E}_k)$ the time-varying graph with $k\in\mathbb{N}$ being the time index. Define the union of graphs $\mathcal{G}_l=(\mathcal{V},\mathcal{E}_l),l=1,\ldots,m$ as $\cup_{l=1}^m\mathcal{G}_l=(\mathcal{V},\cup_{l=1}^m\mathcal{E}_l)$. At each time $k\in\mathbb{N}$, define an adjacency matrix as $A_k=(a_{ij,k})\in\mathbb{R}^{N\times N}$ such that $a_{ij,k}>0$ if $(j,i)\in\mathcal{E}_k$, and $a_{ij,k}=0$ otherwise. It is standard to assume the following assumptions.

\begin{assumption}[$Q$-Strong-Connectivity and Weights Rule]\label{a2}
~
\begin{enumerate}
  \item $\mathcal{G}_k$ is {\em uniformly jointly strongly connected}, i.e., there exists an integer $Q>0$ such that the graph union $\cup_{l=1}^Q \mathcal{G}_{k+l}$ is strongly connected for all $k\geq 0$.
  \item For all $k\in\mathbb{N}$, $A_k$ is doubly stochastic, i.e., $\sum_{j=1}^N a_{ij,k}=1$ and $\sum_{i=1}a_{ij,k}=1$ for all $i,j\in[N]$. Meanwhile, $a_{ij,k}>\underline{a}$ for a constant $\underline{a}\in(0,1)$ whenever $a_{ij,k}>0$, and $a_{ii,k}>\underline{a}$ for all $i\in[N]$, $k\in\mathbb{N}$.
\end{enumerate}
\end{assumption}

%

\section{Main results}\label{s3}

This section aims at proposing two distributed algorithms for coping with problem (\ref{6}), i.e., distributed KM iteration and distributed block-coordinate KM iteration.

\subsection{The D-KM iteration}\label{s3.1}

For problem (\ref{6}), if there exists a global computing unit who knows the exact information $F$, then one famous centralized algorithm, called the KM iteration \cite{cominetti2014rate,matsushita2017convergence}, can be used
\begin{align}
x_{k+1}=x_k+\alpha_k (F(x_k)-x_k),        \label{8}
\end{align}
where $\{\alpha_k\}_{k\in\mathbb{N}}\in[0,1]$ is a sequence of relaxation parameters. It can be found in \cite{reich1979weak} that the KM iteration is weakly convergent to a fixed point of $F$ under mild conditions, e.g., $\sum_{j=1}^\infty \alpha_j(1-\alpha_j)=\infty$. Nevertheless, the classical centralized algorithms are not applicable here since no global coordinators can access the total information on $F$, instead it is imperative to develop distributed algorithms, which are based on local information exchange among all agents.

Inspired by the classical KM iteration, a distributed KM (D-KM) iteration is proposed as
\begin{subequations}
\begin{align}
\hat{x}_{i,k}&=\sum_{j=1}^N a_{ij,k}x_{j,k},          \label{9a}\\
x_{i,k+1}&=\hat{x}_{i,k}+\alpha_{k}(F_i(\hat{x}_{i,k})-\hat{x}_{i,k}),~~i\in[N]              \label{9b}
\end{align}               \label{9}
\end{subequations}
where $\hat{x}_{i,k}$ stands for an aggregate message received from the neighbors of agent $i$ at time instant $k>0$, $x_{i,k}$ is an estimate of a fixed point of the global operator $F$ by agent $i$ at time slot $k$ for all $i\in[N]$, and $\alpha_k$ is a positive stepsize. As for (\ref{9}), it is easy to see that the algorithm is distributed since every agent only makes use of its local information.

To proceed, it is useful to postulate a few properties on $\alpha_k$.
\begin{assumption}[Stepsizes]\label{a3}
The stepsize $\alpha_k$ is nonincreasing for each $k\in\mathbb{N}$, satisfying the following properties:
\begin{align}
\alpha_k\in (0,1],~~~\sum_{k=1}^\infty\alpha_k=\infty,~~~\sum_{k=1}^\infty\alpha_k\alpha_{\lfloor\frac{k}{2}\rfloor}<\infty.          \label{10}
\end{align}
\end{assumption}

It is now ready for us to present the main result on the D-KM iteration (\ref{9}).

\begin{theorem}\label{t1}
If Assumptions \ref{a1}-\ref{a3} hold, then all $x_{i,k}$'s generated by the D-KM iteration (\ref{9}) are bounded and converge weakly to a common point in $Fix(F)$.
\end{theorem}
{\bf Proof.}
The proof is postponed to Appendix \ref{s7.1}.
\hfill\rule{2mm}{2mm}

\begin{remark}\label{rm2}
In contrast to \cite{fullmer2018distributed}, where the common fixed point seeking problem is studied in the Euclidean space, the problem here is more general, where we do not assume a common fixed point for all local operators and consider the problem in real Hilbert spaces. Our work also extends the case where all local operators have nonempty common fixed points addressed in \cite{li2019distributed3} in real Hilbert spaces. Note that the methods employed in \cite{li2019distributed3} are no longer applicable to the problem studied in this paper. Also, without the property of having at least a common fixed point for all local operators, the problem here requires more stringent assumptions than those in \cite{li2019distributed3}. That is, Assumption \ref{a1} and the column-stochasticity in Assumption \ref{a2} in this paper are not needed in \cite{li2019distributed3}, and the conditions on stepsize $\{\alpha_k\}$ in Assumption \ref{a3} are more restrictive that that in \cite{li2019distributed3}, where $\alpha_k$ can even be a positive constant.
\end{remark}

\subsection{The D-BKM iteration}\label{s3.2}

This subsection is concerned with developing another distributed algorithm for handling problem (\ref{6}). In reality, it may be computationally expensive or prohibitive to compute the whole coordinates of $F_i$ at a point for $i\in[N]$. Instead, it is preferable and practical to randomly compute only a partial coordinates of $F_i$ at the point for $i\in[N]$, in order to alleviate the computational complexity. For example, in Remark 1 the operator $T=Id-\tau\nabla f$ involves the gradient $\nabla f$, which is known to be computationally heavy to calculate the entire coordinates of $\nabla f(x)$ at a point $x$ especially when $x$ is of large dimension, while it is easier to only compute some partial coordinates of $\nabla f(x)$.

In this subsection, $\mathcal{H}=\mathcal{H}_1\oplus\cdots\oplus\mathcal{H}_m$ is the direct Hilbert sum with Borel $\sigma$-algebra $\mathcal{B}$ and each $\mathcal{H}_i,i\in[m]$ being a separable real Hilbert space. $\mathcal{H}$ is endowed with the same inner product $\langle\cdot,\cdot\rangle$ and associated norm $\|\cdot\|$ as before. Denote by $x=(x_1,\ldots,x_m)$ a generic vector in $\mathcal{H}$. All local operators are partitioned into $m$ block-coordinates, i.e., $F_i:x\mapsto (F_{i1}(x),\ldots,F_{im}(x))$ with $F_{il}:\mathcal{H}\to\mathcal{H}_l$ being measurable for all $i\in[N]$ and $l\in[m]$.

To tackle problem (\ref{6}), another distributed algorithm, called distributed block-coordinate KM (D-BKM) iteration, is designed in (\ref{01}). Let each initial vector $x_{i,0}$ be a $\mathcal{H}$-valued random variable for all $i\in[N]$. At time $k+1$, a global coordinator randomly selects a block-coordinate number, say $q\in[m]$, with a probability distribution, and then broadcast the number $q$ to all agents. Subsequently, each agent $i\in[N]$ only needs to compute $F_{iq}(\hat{x}_{i,k})$ (instead of the whole coordinates $F_i(\hat{x}_{i,k})$), updates all block-coordinates as in (\ref{01}) with $b_{l,k}=1$ for $l=q$ and with $b_{l,k}=0$ for $l\neq q$, and then propagates $x_{il,k+1},l\in[m]$ to its out-neighbors.
\begin{subequations}
\begin{align}
\hat{x}_{il,k}&=\sum_{j=1}^N a_{ij,k}x_{jl,k},~~~\forall l\in[m],i\in[N]                         \label{01a}\\
x_{il,k+1}&=\hat{x}_{il,k}+b_{l,k}\alpha_{k}(F_{il}(\hat{x}_{i,k})-\hat{x}_{il,k}),     \label{01b}
\end{align}                  \label{01}
\end{subequations}
where $x_{i,k}=(x_{i1,k},\ldots,x_{im,k})$ represents an estimate of a fixed point of the global operator $F$ for agent $i$ at time step $k$, and $\hat{x}_{il,k}$ is the $l$-th block-coordinate of $\hat{x}_{i,k}=(\hat{x}_{i1,k},\ldots,\hat{x}_{im,k})$, serving as an aggregate information collected from its neighbors at time instant $k$. Moreover, $\{b_{l,k}\}_{k\in\mathbb{N}}$ is a sequence of $\{0,1\}$-valued random variables (independently identically distributed). Additionally, $\{\alpha_{k}\}_{k\in\mathbb{N}}$ is a sequence of stepsizes.

For brevity, let $\chi_{i,k}:=\sigma(x_{i,0},\ldots,x_{i,k})$, $\chi_k:=\sigma(\chi_{1,k},\ldots,\chi_{N,k})$, and $\mathcal{E}_{i,k}:=\sigma(b_{i,k})$ for $i\in[N]$ and $k\in\mathbb{N}$. Set $\chi=\{\chi_k\}_{k\in\mathbb{N}}$. It is standard to assume that $\mathcal{E}_{i,k}$ is independent of both $\chi_k$ and $\mathcal{E}_{j,k}$ for $j\neq i\in[N]$. Moreover, define $p_l:=\mathbb{P}(b_{l,0}=1)$ for $l\in[m]$, and assume $p_l>0$ for all $l\in[m]$, connoting that there is an opportunity for every block-coordinate to be selected.

To proceed, it is helpful for the following analysis to rewrite (\ref{01}) as
\begin{align}
x_{il,k+1}=\hat{x}_{il,k}+\alpha_{k}(T_{il,k}-\hat{x}_{il,k}),      \label{02}
\end{align}
where
\begin{align}
T_{il,k}:=\hat{x}_{il,k}+b_{l,k}(F_{il}(\hat{x}_{i,k})-\hat{x}_{il,k}).            \label{03}
\end{align}

By defining $T_{i,k}:=(T_{i1,k},\ldots,T_{im,k})$, (\ref{02}) can be rewritten in a compact form as
\begin{align}
x_{i,k+1}=\hat{x}_{i,k}+\alpha_{k}(T_{i,k}-\hat{x}_{i,k}).      \label{04}
\end{align}

Summing (\ref{04}) over $i\in[N]$ and invoking the double-stochasticity of $A_k$ in Assumption \ref{a2}, it can be obtained that
\begin{align}
\bar{x}_{k+1}=\bar{x}_{k}+\alpha_{k}\Big(\frac{\sum_{i=1}^N T_{i,k}}{N}-\bar{x}_{k}\Big),      \label{05a}
\end{align}
where $\bar{x}_{k}:=\sum_{i=1}^N x_{i,k}/N$.

Before presenting the main result, it is helpful to propose a new norm $|||\cdot|||$ and corresponding inner product $\langle\langle\cdot,\cdot\rangle\rangle$ on $\mathcal{H}$ as in \cite{li2019distributed3}, i.e., for any $y,z\in\mathcal{H}$
\begin{align}
|||y|||^2:=\sum_{l=1}^m \frac{1}{p_l}\|y_l\|^2,~~~\langle\langle y,z\rangle\rangle :=\sum_{l=1}^m \frac{1}{p_l}\langle y_l,z_l\rangle,              \label{09}
\end{align}
for which it is easy to verify that $\|y\|^2\leq |||y|||^2\leq \|y\|^2/p_0$, where $p_0:=\min_{l\in[m]}p_l$.

With the above at hand, we are now in a position to give the main result on (\ref{01}).

\begin{theorem}\label{t2}
If Assumptions \ref{a1}-\ref{a3} hold, then all $x_{i,k}$'s generated by the D-BKM iteration (\ref{01}) are bounded and converge weakly, in the space $(\mathcal{H},|||\cdot|||)$, to a common point in $Fix(F)$ a.s.
\end{theorem}
{\bf Proof.}
The proof is postponed to Appendix \ref{s7.2}.
\hfill\rule{2mm}{2mm}

\begin{remark}\label{rm3}
It is worthwhile to notice that the D-BKM iteration (\ref{01}) requires a global coordinator to randomly select one of $m$ block-coordinates for all agents at each time instant, while no such sort of global coordinator is required in \cite{li2019distributed3} for the case where all local operators have common fixed points. In this regard, the main difficulty lies in the general setup in this paper that no common fixed point is assumed for all local operators, and it is still open to consider the case where no such coordinator is present in the network, which is left as one of our future works.
\end{remark}

\section{Applications}\label{s4}

This section focuses on applying the previous results to several concrete problems, i.e., distributed optimization, a distributed shortest distance problem, and the problem of solving linear algebraic equations in a distributed manner.

\subsection{Distributed gradient descent algorithm}\label{s4.1}

In this subsection, let us consider the following distributed optimization problem in real Hilbert spaces
\begin{align}
\text{min}~~~f(x):=\sum_{i=1}^N f_i(x),         \label{x1}
\end{align}
where $f_i:\mathcal{H}\to\mathbb{R}$ is a convex and differentiable function with Lipschitz gradient of constant $L_i$ for all $i\in[N]$.

Define $L:=\max_{i\in[N]}L_i$ and $F_i:=Id-\tau\nabla f_i$ for a constant $\tau\in(0,2/L)$. It can be then asserted that the operators $F_i$'s are all nonexpansive \cite{liang2016convergence}. As a result, it is easy to verify that problem (\ref{x1}) is equivalent to finding fixed points of the operator $F:=\sum_{i=1}^N F_i/N$, which is exactly the same as (\ref{6}).

In this respect, the D-KM and D-BKM iterations become
\begin{align}
x_{i,k+1}&=\hat{x}_{i,k}+\beta_k \nabla f_i(\hat{x}_{i,k}),                  \label{x2}\\
x_{il,k+1}&=\hat{x}_{il,k}+b_{l,k}\beta_k \nabla f_{il}(\hat{x}_{i,k}),         \label{x3}
\end{align}
respectively, of which (\ref{x2}) is the classical distributed gradient descent (DGD) algorithm \cite{nedic2009distributed,liu2017convergence}, and (\ref{x3}) is a block-coordinate based DGD algorithm. Wherein, $\beta_k:=\tau\alpha_k$, $\alpha_k$ is the stepsize satisfying Assumption \ref{a3}, $\hat{x}_{i,k}$ and $b_{l,k}$ are the same as in (\ref{9a}) and (\ref{01b}), respectively, and $f_{il},\hat{x}_{il,k}$ are the $l$-th components of $f_i$ and $\hat{x}_{i,k}$ for $l\in[m]$, respectively. Hence, under Assumptions \ref{a1}-\ref{a3}, Theorems \ref{t1} and \ref{t2} hold for (\ref{x2}) and (\ref{x3}), respectively. It should be noted that Assumption \ref{a1} in this case amounts to saying that the gradients of all $f_i$'s are bounded, which is widely employed in the literature, e.g., \cite{nedic2009distributed,liu2017convergence,mateos2017distributed,xie2018distributed,liang2019distributed}.

\begin{remark}\label{rm4}
To our best knowledge, (\ref{x3}) is the first block-coordinate DGD algorithm for distributed optimization under time-varying directed communication graphs. We note that the distributed optimization (\ref{x1}) in real Hilbert spaces was also studied recently in \cite{alaviani2019distributed}. However, problem (\ref{x1}) is a special case of problem (\ref{6}), which is a more general problem.
\end{remark}

\subsection{A distributed shortest distance problem}\label{s4.2}

Consider now a special case of (\ref{6}) by selecting $F_i=P_{X_i}$, where $X_i$ is a subset of $\mathcal{H}$ for all $i\in[N]$, and all $X_i$'s are convex and compact. In this case, problem (\ref{6}) degenerates to
\begin{align}
\text{find}~x\in\mathcal{H},~~\text{s.t.}~~x\in Fix(F),~~F:=\frac{1}{N}\sum_{i=1}^N P_{X_i},            \label{x4}
\end{align}
which is called a {\em distributed shortest distance} problem, since (\ref{x4}) is equivalent to minimizing $\sum_{i=1}^N d_{X_i}^2(x)$ for $x\in\mathcal{H}$ that can be applied such as in source localization \cite{zhang2015distributed}.

Along this line, the D-KM and D-BKM iterations become
\begin{align}
x_{i,k+1}&=\hat{x}_{i,k}+\alpha_k (P_{X_i}(\hat{x}_{i,k})-\hat{x}_{i,k}),              \label{x6}\\
x_{il,k+1}&=\hat{x}_{il,k}+b_{l,k}\alpha_k (P_{X_il}(\hat{x}_{i,k})-\hat{x}_{il,k}),         \label{x7}
\end{align}
respectively, where $P_{X_il}$ is the $l$-th component of $P_{X_i}$, i.e., $P_{X_i}=(P_{X_i1},\ldots,P_{X_im})$.

Due to the compactness of $X_i$'s, it is straightforward to see that $\|P_{X_i}(x)\|$ is bounded for any $x\in\mathcal{H}$ and all $i\in[N]$, which however cannot ensure the correctness of Assumption \ref{a1}. Thus, there is a need to first show the boundedness in Assumption \ref{a1}.

\begin{proposition}\label{p1}
Under Assumptions \ref{a2} and \ref{a3}, for iterations (\ref{x6}) and (\ref{x7}) with $p_l=1/m$ for all $l\in[m]$, Assumption \ref{a1} holds.
\end{proposition}
{\bf Proof.}
The proof is postponed to Appendix \ref{s7.3}.
\hfill\rule{2mm}{2mm}

Equipped with Proposition \ref{p1}, Assumptions \ref{a2} and \ref{a3}, the weak convergence to a solution of (\ref{x4}) for iterations (\ref{x6}) and (\ref{x7}) with $p_l=1/m$ can be ensured by Theorems \ref{t1} and \ref{t2}.

\begin{remark}\label{rm5}
The distributed shortest distance problem has been investigated in \cite{lou2016distributed,lin2018multiagent}, where distributed continuous-time algorithms are designed in the Euclidean space. In contrast, this paper develops distributed discrete-time algorithms (\ref{x6}) and (\ref{x7}) in the Hilbert space. Particularly, to our best knowledge, (\ref{x7}) is the first distributed block-coordinate algorithm for (\ref{x4}) under time-varying directed graphs. It is also worth mentioning that problem (\ref{x4}) encompasses the classical {\em convex feasibility} problem \cite{necoara2018randomized,kruger2018set} as a special case, i.e., $\text{find}~x\in\mathcal{H},~~\text{s.t.}~~x\in \cap_{i=1}^N X_i$, when all $X_i$'s have a common point due to $Fix(F)=\cap_{i=1}^N Fix(F_i)=\cap_{i=1}^N X_i$ (e.g., Proposition 4.47 in \cite{bauschke2017convex}).
\end{remark}

\subsection{Solving linear algebraic equations}\label{s4.3}

This subsection is on a classical problem of solving a linear algebraic equation in a distributed approach, that is,
\begin{align}
\text{solve}~~~H(x):=Rx-r=0,           \label{x8}
\end{align}
where $x\in\mathbb{R}^n$ is the decision vector, $r\in\mathbb{R}^n$, and $R\in\mathbb{R}^{n\times n}$ is a symmetric positive semi-definite matrix. This setting is general, since the form $Ax=b$ with $A\in\mathbb{R}^{d\times n},b\in\mathbb{R}^d$ can be cast as (\ref{x8}) by defining $R=A^\top A$ and $r=A^\top b$.

It is easy to see that (\ref{x8}) amounts to finding a fixed point of the operator $F$, where $F$ is defined by
\begin{align}
F: x\mapsto (I-\theta R)x+\theta r,~~~\forall x\in\mathbb{R}^n       \label{x9}
\end{align}
for any constant $\theta>0$.

Consider the case when $R$ and $r$ are separable, i.e., there are symmetric positive semi-definite matrices $R_i$'s and vectors $r_i$'s such that
\begin{align}
R=\frac{\sum_{i=1}^N R_i}{N},~~~r=\frac{\sum_{i=1}^N r_i}{N},            \label{x10}
\end{align}
where $\theta\in(0,2/\lambda_M]$, $\lambda_M:=\max_{i\in[N]}\{\lambda_{max}(R_i)\}$, and $\lambda_{max}(\cdot)$ means the largest eigenvalue of a matrix. In this case, agent $i$ is only capable of accessing the information on $R_i$ and $r_i$ for each $i\in[N]$.

Consequently, the global operator $F$ can be rewritten as
\begin{align}
F=\frac{\sum_{i=1}^N F_i}{N},           \label{x11}
\end{align}
where each local operator $F_i$ is defined as
\begin{align}
F_i: x\mapsto (I-\theta R_i)x+\theta r_i,~~~\forall x\in\mathbb{R}^n           \label{x12}
\end{align}
which is nonexpansive due to $\theta\in(0,2/\lambda_M]$.

At this step, it is easy to observe that problem (\ref{x8}) is eventually equivalent to
\begin{align}
\text{find}~~~x\in\mathbb{R}^n,~~s.t.~~x\in Fix(F),~~F=\frac{\sum_{i=1}^N F_i}{N},           \label{x13}
\end{align}
which fits exactly problem (\ref{6}).

However, Assumption \ref{a1} cannot be directly confirmed for $F_i$'s in (\ref{x12}), for which the following result indicates that Assumption \ref{a1} is indeed correct for the D-KM iteration.
\begin{proposition}\label{p2}
Under Assumptions \ref{a2} and \ref{a3}, for the D-KM iteration with $F_i$'s in (\ref{x12}), Assumption \ref{a1} holds.
\end{proposition}
{\bf Proof.}
The proof is postponed to Appendix \ref{s7.4}.
\hfill\rule{2mm}{2mm}

Based on Proposition \ref{p2}, Assumptions \ref{a2} and \ref{a3}, Theorem \ref{t1} holds, that is, the D-KM iteration is weakly convergent to a solution of problem (\ref{x8}).

\begin{remark}\label{rm6}
Lots of existing works focus on solving (\ref{x8}) in a distributed manner, such as \cite{mou2015distributed,wang2017further,wang2019distributed,alaviani2018distributed}. However, these works study the case where each agent only knows some complete rows of $R$ and corresponding entries of $r$, while each agent here can know a sub-matrix of $R$ as given in (\ref{x10}), thus providing a new perspective for solving (\ref{x8}) in a distributed manner. Although a relatively general partition of $R$ is addressed in \cite{wang2019scalable} based on a double-layered network, it only applies to undirected and fixed graphs, and only provides continuous-time algorithms. In contrast, the D-KM iteration here is a discrete-time algorithm, and the underlying communication graph is simple one-layered, directed, and time-varying.
\end{remark}

\section{Numerical examples}\label{s5}

In this section, two numerical examples are given to corroborate the developed D-KM and D-BKM iterations in this paper. To this end, let us focus on the distributed shortest distance problem (\ref{x4}) in Section \ref{s4.2}, i.e., setting $F_i=P_{X_i}$ for all $i\in[N]$, where $X_i$'s are convex and compact subsets of $\mathbb{R}^3$ and each $X_i$ is only accessible to agent $i$. Let $\alpha_k=1/k^{0.7}$ and $X_i=[\sqrt{i},\sqrt{i+1}]\times [\sin(i\pi/2),1+\sin(i\pi/2)]\times [\sqrt{i}-\sqrt{N}+2,\sqrt{i}]$ for all $i\in[N]$.

\begin{figure}[H]
\centering
\subfigure[]{\includegraphics[width=0.8in]{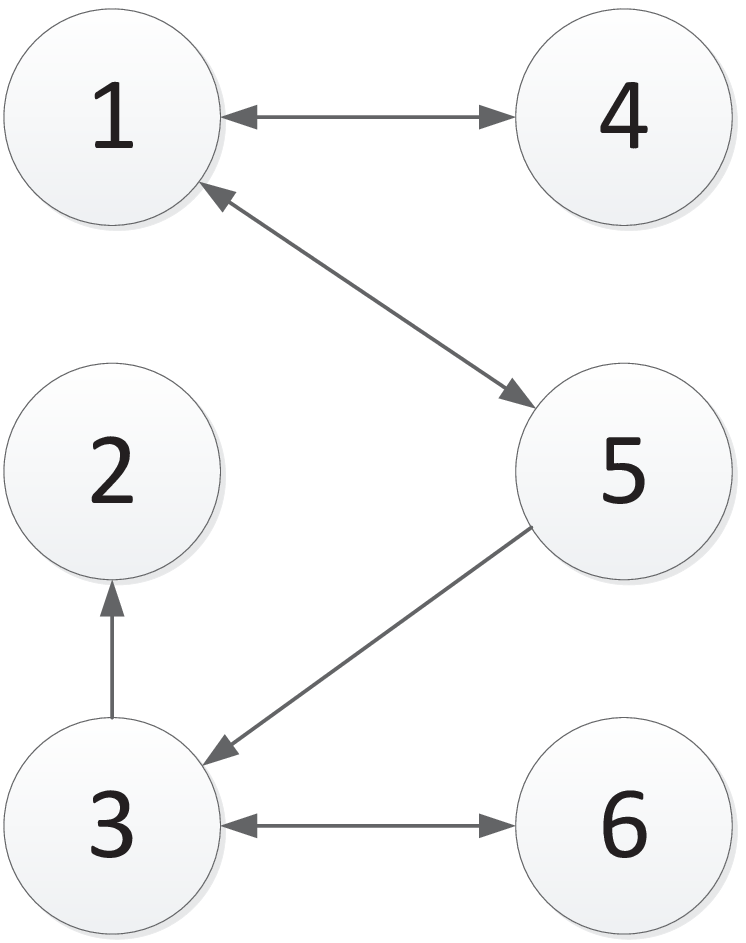}}\hspace{1.1cm}
\subfigure[]{\includegraphics[width=0.8in]{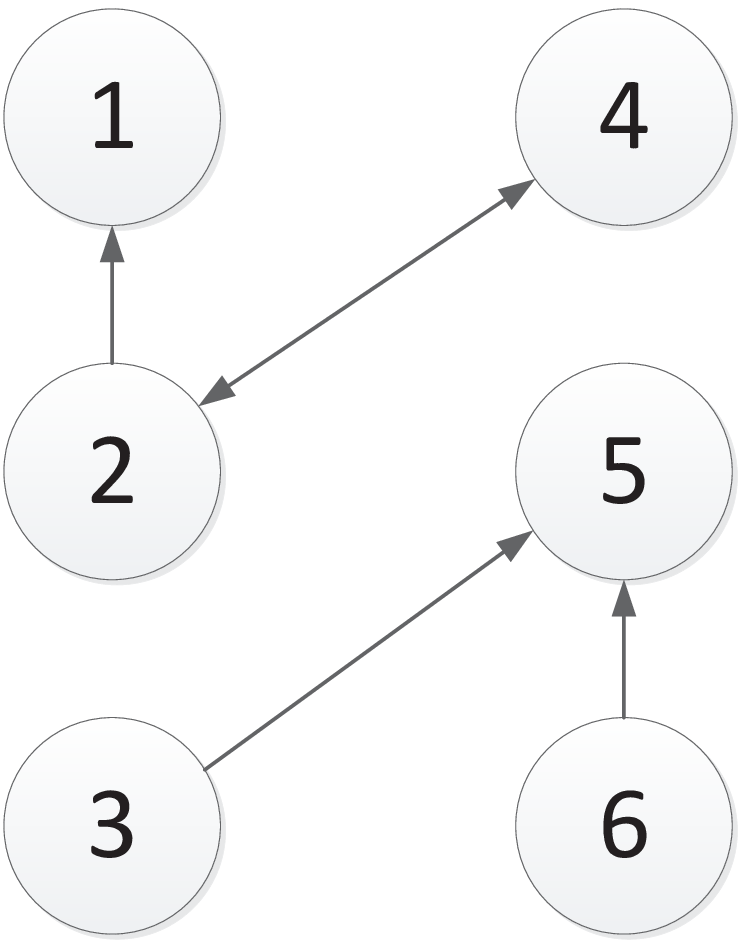}}
\caption{Two switching graphs for the D-KM iteration.}
\label{f1}
\end{figure}

\begin{figure}[H]
\centering
\includegraphics[width=3.1in]{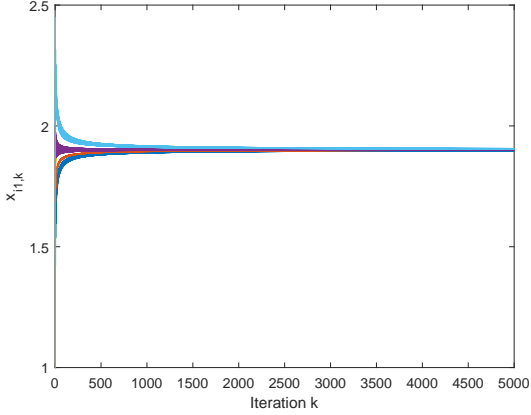}
\caption{Evolutions of $x_{i1,k}$'s by the D-KM iteration.}
\label{f2}
\end{figure}

\begin{figure}[H]
\centering
\includegraphics[width=3.1in]{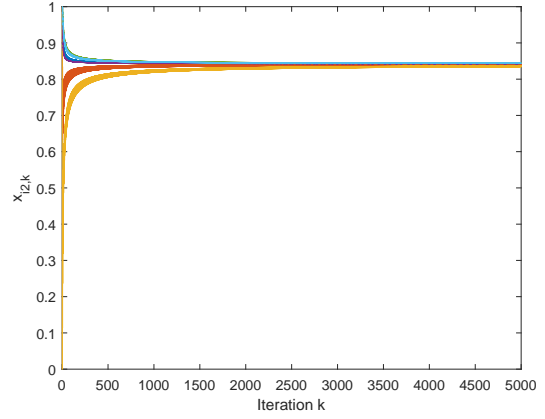}
\caption{Evolutions of $x_{i2,k}$'s by the D-KM iteration.}
\label{f3}
\end{figure}

\begin{figure}[H]
\centering
\includegraphics[width=3.1in]{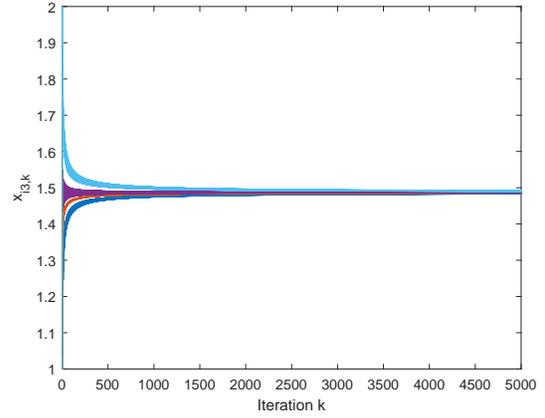}
\caption{Evolutions of $x_{i3,k}$'s by the D-KM iteration.}
\label{f4}
\end{figure}

With regard to the D-KM iteration, consider $N=6$ agents in the network. The communication graphs is uniformly jointly strongly connected with $Q=2$, for which two switching graphs are shown in Fig. \ref{f1}. In this setup, by choosing the initial vectors arbitrarily, the simulation trajectories are given in Figs. \ref{f2}-\ref{f4} for the first, second, and third coordinates, respectively, indicating that the D-KM iteration is indeed convergent to a solution of the studied problem as asserted in Theorem \ref{t1}.

\begin{figure}[H]
\centering
\includegraphics[width=3.1in]{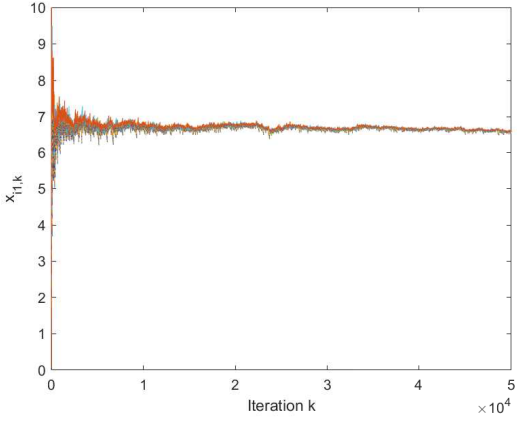}
\caption{Evolutions of $x_{i1,k}$'s by the D-BKM iteration.}
\label{f5}
\end{figure}

\begin{figure}[H]
\centering
\includegraphics[width=3.1in]{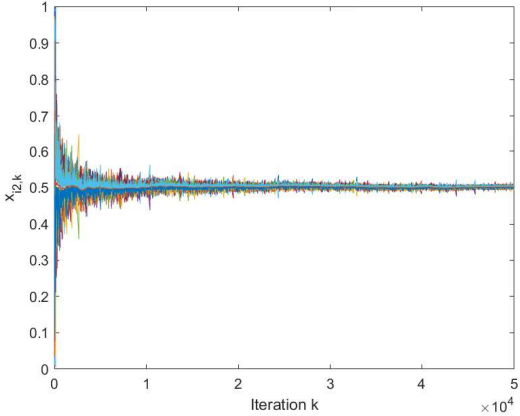}
\caption{Evolutions of $x_{i2,k}$'s by the D-BKM iteration.}
\label{f6}
\end{figure}

\begin{figure}[H]
\centering
\includegraphics[width=3.1in]{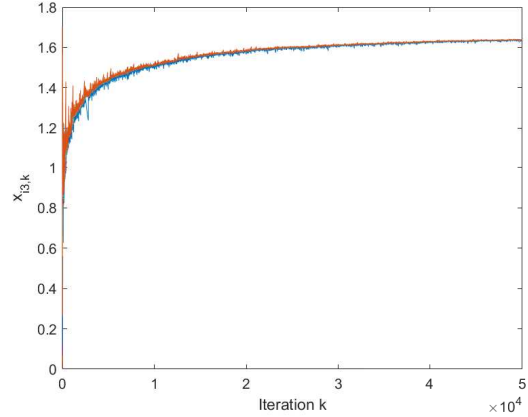}
\caption{Evolutions of $x_{i3,k}$'s by the D-BKM iteration.}
\label{f7}
\end{figure}

For the D-BKM iteration, let us take into account a network of $N=100$ agents and the connectivity parameter $Q=10$. Set $p_l=\mathbb{P}(b_{l,0}=1)=1/3$ for $l\in[3]$. For any randomly selected initial vectors, performing the D-BKM iteration yields the evolutions of all $x_{i,k}$'s in Figs. \ref{f5}-\ref{f7} for each coordinate. In this scenario, the simulation results still support the theoretical claim in Theorem \ref{t2}.

\section{Conclusion}\label{s6}

This paper has addressed the fixed point seeking problem for a global nonexpansive operator in real Hilbert spaces, where the global operator is a sum of local nonexpansive operators and each local operator is only privately accessible to individual agent. In the setup, all agents must communicate via information exchange to handle the problem in a cooperative way, and the communications among all agents are directed and time-varying, satisfying $Q$-strong-connectivity. To deal with the problem, two distributed algorithms have been developed and rigorously analyzed for the weak convergence. One algorithm is called the D-KM iteration, motivated by the classical centralized KM iteration, and the other is called the D-BKM iteration, which is the block-coordinate version of the D-KM iteration. As applications of the theoretical results in this paper, three problems have been considered, i.e., distributed optimization, a distributed shortest distance problem, and solving linear algebraic equations in a distributed manner. Numerical examples have also been presented to support the theoretical results. Directions of future work can be put on the convergence speed, fully distributed block-coordinate KM iteration (i.e., without the global coordinator for selecting the updated block-coordinate for all agents), and the case with random communication graphs.

%

\section{Appendix}\label{s7}

\subsection{Proof of Theorem \ref{t1}}\label{s7.1}

To begin, it is useful to list some lemmas.

\begin{lemma}[\cite{li2019distributed3}]\label{l1}
Let $\{v_k\}$ be a sequence of nonnegative scalars such that for all $k\geq 0$
\begin{align}
v_{k+1}\leq(1+b_k)v_k-u_k+c_k,      \nonumber
\end{align}
where $b_k\geq0$, $u_k\geq0$ and $c_k\geq0$ for all $k\geq0$ with $\sum_{k=1}^\infty b_k<\infty$ and $\sum_{k=1}^\infty c_k<\infty$. Then, the sequence $\{v_k\}$ converges to some $v\geq 0$ and $\sum_{k=1}^\infty u_k<\infty$.
\end{lemma}

\begin{lemma}[\cite{bauschke2017convex}]\label{l2}
Let $x,y\in\mathcal{H}$, and let $r\in\mathbb{R}$. Then
\begin{align}
\|rx+(1-r)y\|^2&=r\|x\|^2+(1-r)\|y\|^2      \nonumber\\
&\hspace{0.4cm}-r(1-r)\|x-y\|^2.       \nonumber
\end{align}
\end{lemma}

\begin{lemma}[\cite{li2019distributed3}]\label{l3}
Let $A\in\mathbb{R}^{n\times n}$ and $B$ be a linear operator in real Hilbert space $\mathcal{H}$, then $\|A\otimes B\|\leq na_{max}\|B\|$, where $a_{max}$ is the largest entry of the matrix $A$ in the modulus sense.
\end{lemma}

\begin{lemma}\label{l4}
Under Assumptions \ref{a1}-\ref{a3}, there holds that
\begin{align}
\|x_{i,k}-\bar{x}_k\|=O(\alpha_{\lfloor\frac{k}{2}\rfloor}),        \label{pf2}
\end{align}
where $\bar{x}_k=\sum_{i=1}^N x_{i,k}/N$ for all $k\in\mathbb{N}$.
\end{lemma}
{\bf Proof.}
It is easy to see that the D-KM iteration (\ref{9}) can be compactly written as
\begin{align}
x_{k+1}=(A_k\otimes Id)x_k+\epsilon_k,          \label{pf1}
\end{align}
where $x_k:=col(x_{1,k},\ldots,x_{N,k})$ and $\epsilon_k:=col(\alpha_k(F_1(\hat{x}_{1,k})-\hat{x}_{1,k}),\ldots,\alpha_k(F_N(\hat{x}_{N,k})-\hat{x}_{N,k}))$. Note that $\|\epsilon_k\|\to 0$ as $\alpha_k\to 0$ and $F_i(\hat{x}_{i,k})-\hat{x}_{i,k}$ are bounded under Assumption \ref{a1}. For (\ref{pf1}), invoking the same reasoning as in Lemmas 3 and 4 in \cite{xie2018distributed}, together with Lemma \ref{l3}, yields the conclusion.
\hfill\rule{2mm}{2mm}

Equipped with the above lemmas, it is now ready to prove Theorem \ref{t1}.

{\em Proof of Theorem \ref{t1}:}
For (\ref{pf1}), by left multiplying a row vector $(1,\ldots,1)$, it is easy to obtain that
\begin{align}
\bar{x}_{k+1}=\bar{x}_k+\alpha_k \Big(\frac{\sum_{i=1}^N F_i(\hat{x}_{i,k})}{N}-\bar{x}_k\Big).       \label{pf3}
\end{align}
Now, for any $x^*\in Fix(F)$, appealing to (\ref{pf3}) gives rise to
\begin{align}
&\|\bar{x}_{k+1}-x^*\|           \nonumber\\
&=\|(1-\alpha_k)(\bar{x}_k-x^*)+\alpha_k \Big(\frac{\sum_{i=1}^N F_i(\hat{x}_{i,k})}{N}-x^*\Big)\|         \nonumber\\
&\leq (1-\alpha_k)\|\bar{x}_k-x^*\|+\alpha_k\|\frac{\sum_{i=1}^N F_i(\bar{x}_{k})}{N}-x^*\|             \nonumber\\
&\hspace{0.4cm}+\alpha_k\|\frac{\sum_{i=1}^N F_i(\hat{x}_{i,k})}{N}-\frac{\sum_{i=1}^N F_i(\bar{x}_{k})}{N}\|           \nonumber\\
&\leq \|\bar{x}_k-x^*\|+\frac{\alpha_k}{N}\sum_{i=1}^N\|F_i(\hat{x}_{i,k})-F_i(\bar{x}_{k})\|               \nonumber\\
&\leq \|\bar{x}_k-x^*\|+\frac{\alpha_k}{N}\sum_{i=1}^N\|\hat{x}_{i,k}-\bar{x}_{k}\|,                   \label{pf4}
\end{align}
where the nonexpansiveness of $F_i$'s have been utilized in the last two inequalities.

Let us consider the term $\|\hat{x}_{i,k}-\bar{x}_{k}\|$. One has that
\begin{align}
\|\hat{x}_{i,k}-\bar{x}_{k}\|&=\|\sum_{j=1}^N a_{ij,k}x_{j,k}-\bar{x}_{k}\|        \nonumber\\
&\leq \sum_{j=1}^N a_{ij,k}\|x_{j,k}-\bar{x}_k\|.                                \label{pf5}
\end{align}
In view of $\sum_{i=1}^N a_{ij,k}=1$, inserting (\ref{pf5}) to (\ref{pf4}) leads to
\begin{align}
\|\bar{x}_{k+1}-x^*\|&\leq \|\bar{x}_k-x^*\|+\frac{\alpha_k}{N}\sum_{j=1}^N\|\hat{x}_{j,k}-\bar{x}_{k}\|                \nonumber\\
&\leq \|\bar{x}_k-x^*\|+c_1\alpha_k\alpha_{\lfloor\frac{k}{2}\rfloor},                             \label{pf6}
\end{align}
where $c_1>0$ is some constant, and the last inequality has employed (\ref{pf2}).

By Assumption \ref{a3}, applying Lemma \ref{l1} to (\ref{pf6}) yields that
\begin{align}
\|\bar{x}_k-x^*\|~\text{is convergent}                  \label{pf7}
\end{align}
and hence $\bar{x}_k,x_{i,k}$ are bounded for all $k\in\mathbb{N}$.

In the following, it remains to prove the weak convergence of $x_{i,k}$'s. To do so, one can obtain that
\begin{align}
&\|\bar{x}_{k+1}-x^*\|^2     \nonumber\\
&=\|(1-\alpha_k)(\bar{x}_k-x^*)+\alpha_k\Big(\frac{\sum_{i=1}^N F_i(\hat{x}_{i,k})}{N}-x^*\Big)\|^2,        \nonumber
\end{align}
which, invoking Lemma \ref{l2}, follows that
\begin{align}
&\|\bar{x}_{k+1}-x^*\|^2     \nonumber\\
&=(1-\alpha_k)\|\bar{x}_k-x^*\|^2+\alpha_k\underbrace{\Big\|\frac{\sum_{i=1}^N F_i(\hat{x}_{i,k})}{N}-x^*\Big\|^2}_{=:S_1}      \nonumber\\
&\hspace{0.4cm}-\alpha_k(1-\alpha_k)\underbrace{\Big\|\frac{\sum_{i=1}^N F_i(\hat{x}_{i,k})}{N}-\bar{x}_k\Big\|^2}_{=:S_2}.        \label{pf8}
\end{align}
Consider terms $S_1$ and $S_2$. For $S_1$, it has that
\begin{align}
S_1&=\Big\|\underbrace{\frac{\sum_{i=1}^N F_i(\bar{x}_k)}{N}-x^*}_{=:a_1}+\underbrace{\frac{\sum_{i=1}^N [F_i(\hat{x}_{i,k})-F_i(\bar{x}_k)]}{N}}_{=:a_2}\Big\|^2     \nonumber\\
&=\|a_1\|^2+\|a_2\|^2+2\langle a_1,a_2\rangle           \nonumber\\
&\leq \|a_1\|^2+\|a_2\|^2+2\|a_1\|\cdot\|a_2\|           \nonumber\\
&\leq \|\bar{x}_k-x^*\|^2+\frac{2\|a_1\|}{N}\sum_{i=1}^N \|\hat{x}_{i,k}-\bar{x}_k\|           \nonumber\\
&\hspace{0.4cm}+\frac{1}{N}\sum_{i=1}^N \|\hat{x}_{i,k}-\bar{x}_k\|^2       \nonumber\\
&\leq \|\bar{x}_k-x^*\|^2+\frac{2\|a_1\|}{N}\sum_{j=1}^N \|x_{j,k}-\bar{x}_k\|          \nonumber\\
&\hspace{0.4cm}+\frac{1}{N}\sum_{j=1}^N \|x_{j,k}-\bar{x}_k\|^2,               \label{pf9}
\end{align}
where the nonexpansiveness of $F_i$'s are used to get the second inequality, and (\ref{pf5}) is employed to deduce the last inequality.

For $S_2$, it can be obtained that
\begin{align}
S_2&=\Big\|\underbrace{F(\bar{x}_k)-\bar{x}_k}_{=:a_3}+a_2\Big\|^2\geq (\|a_3\|-\|a_2\|)^2          \nonumber\\
&\geq \frac{\|a_3\|^2}{2}-\|a_2\|^2,         \label{pf10}
\end{align}
where the last inequality has applied the fact that $(a-b)^2\geq a^2/2-b^2$ for any $a,b\geq 0$.

Plugging (\ref{pf9}) and (\ref{pf10}) into (\ref{pf8}) results in
\begin{align}
&\|\bar{x}_{k+1}-x^*\|^2     \nonumber\\
&\leq \|\bar{x}_k-x^*\|^2+\frac{2\alpha_k\|a_1\|}{N}\sum_{j=1}^N\|x_{j,k}-\bar{x}_k\|          \nonumber\\
&\hspace{0.4cm}+\frac{\alpha_k}{N}\sum_{j=1}^N \|x_{j,k}-\bar{x}_k\|^2+\frac{\alpha_k(1-\alpha_k)}{N}\sum_{j=1}^N \|x_{j,k}-\bar{x}_k\|^2           \nonumber\\
&\hspace{0.4cm}-\frac{\alpha_k(1-\alpha_k)}{2}\|F(\bar{x}_k)-\bar{x}_k\|^2     \nonumber\\
&\leq \|\bar{x}_k-x^*\|^2-\frac{\alpha_k(1-\alpha_k)}{2}\|F(\bar{x}_k)-\bar{x}_k\|^2       \nonumber\\
&\hspace{0.4cm}+c_2\alpha_k\alpha_{\lfloor\frac{k}{2}\rfloor}        \label{pf11}
\end{align}
for some constant $c_2>0$, where the last inequality has employed (\ref{pf2}) and the boundedness of $\|a_1\|$ due to (\ref{pf7}).

From (\ref{pf11}), invoking Assumption \ref{a3} and Lemma \ref{l1} yields
\begin{align}
\sum_{k=1}^\infty\alpha_k\|F(\bar{x}_k)-\bar{x}_k\|^2<\infty,          \nonumber
\end{align}
which, together with $\sum_{k=1}^\infty \alpha_k=\infty$ in Assumption \ref{a3}, follows
\begin{align}
\mathop{\lim\inf}_{k\to\infty}\|F(\bar{x}_k)-\bar{x}_k\|=0.          \label{pf12}
\end{align}

On the other hand, it can be obtained that
\begin{align}
&\|F(\bar{x}_{k+1})-\bar{x}_{k+1}\|        \nonumber\\
&=\|F(\bar{x}_{k+1})-F(\bar{x}_{k})+F(\bar{x}_{k})-\bar{x}_{k+1}\|           \nonumber\\
&=\|F(\bar{x}_{k+1})-F(\bar{x}_{k})+(1-\alpha_k)(F(\bar{x}_{k})-\bar{x}_k)        \nonumber\\
&\hspace{0.4cm}-\frac{\alpha_k}{N}\sum_{i=1}^N(F_i(\hat{x}_{i,k})-F_i(\bar{x}_k))\|      \nonumber\\
&\leq \|\bar{x}_{k+1}-\bar{x}_k\|+(1-\alpha_k)\|F(\bar{x}_{k})-\bar{x}_k\|              \nonumber\\
&\hspace{0.4cm}+\frac{\alpha_k}{N}\sum_{i=1}^N\|F_i(\hat{x}_{i,k})-F_i(\bar{x}_k)\|      \nonumber\\
&\leq \alpha_k\Big\|\frac{\sum_{i=1}^N F_i(\hat{x}_{i,k})}{N}-\bar{x}_k\Big\|+(1-\alpha_k)\|F(\bar{x}_{k})-\bar{x}_k\|              \nonumber\\
&\hspace{0.4cm}+\frac{\alpha_k}{N}\sum_{i=1}^N\|F_i(\hat{x}_{i,k})-F_i(\bar{x}_k)\|      \nonumber\\
&\leq \|F(\bar{x}_k)-\bar{x}_k\|+\frac{2\alpha_k}{N}\sum_{i=1}^N\|F_i(\hat{x}_{i,k})-F_i(\bar{x}_k)\|         \nonumber\\
&\leq \|F(\bar{x}_k)-\bar{x}_k\|+\frac{2\alpha_k}{N}\sum_{j=1}^N\|x_{j,k}-\bar{x}_k\|         \nonumber\\
&\leq \|F(\bar{x}_k)-\bar{x}_k\|+c_3\alpha_k\alpha_{\lfloor\frac{k}{2}\rfloor},              \label{pf13}
\end{align}
where $c_3>0$ is some constant, (\ref{pf3}) has been leveraged in deriving the second equality and inequality, the nonexpansiveness of $F$ and $F_i$'s have been used in the first and fourth inequalities, and (\ref{pf2}) has been applied to obtain the last inequality.

Applying Lemma \ref{l1} to (\ref{pf13}) implies that $\|F(\bar{x}_k)-\bar{x}_k\|$ is convergent, which in combination with (\ref{pf12}) leads to
\begin{align}
\lim_{k\to\infty}\|F(\bar{x}_k)-\bar{x}_k\|=0.          \label{pf14}
\end{align}
At this moment, for any weak sequential cluster point $x_c$ of $\{\bar{x}_k\}_{k\in\mathbb{N}}$, i.e., $\bar{x}_k\rightharpoonup x_c$, with reference to (\ref{pf14}) and Corollary 4.28 in \cite{bauschke2017convex}, one can claim that $x_c\in Fix(F)$. Subsequently, appealing to Lemma 2.47 in \cite{bauschke2017convex} and (\ref{pf7}) yields that $\bar{x}_k\rightharpoonup x'$ for some point $x'\in Fix(F)$. Therefore, for any $x\in\mathcal{H}$ and $i\in[N]$, one has that
\begin{align}
\langle x_{i,k}-x',x\rangle&=\langle x_{i,k}-\bar{x}_k,x\rangle+\langle \bar{x}_{k}-x',x\rangle            \nonumber\\
&\leq \|x_{i,k}-\bar{x}_k\|\cdot\|x\|+\langle \bar{x}_{k}-x',x\rangle          \nonumber\\
&\to 0.                   \label{pf15}
\end{align}
As a consequence, it can be concluded that $x_{i,k}\rightharpoonup x'$ for all $i\in[N]$. This ends the proof.
\hfill\rule{2mm}{2mm}

\subsection{Proof of Theorem \ref{t2}}\label{s7.2}

The following lemmas are useful for the upcoming analysis.

\begin{lemma}[\cite{robbins1971martin}]\label{l5}
Let $\mathfrak{F}=\{\mathcal{F}_k\}_{k\in\mathbb{N}}$ be a filtration. If $\{z_k\}_{k\in\mathbb{N}}\in\ell_+(\mathfrak{F})$, $\{\varsigma_k\}_{k\in\mathbb{N}}\in\ell_+^1(\mathfrak{F})$, $\{\vartheta_k\}_{k\in\mathbb{N}}\in\ell_+(\mathfrak{F})$, and $\{\eta_k\}_{k\in\mathbb{N}}\in\ell_+^1(\mathfrak{F})$ satisfy the following inequality a.s.:
\begin{align}
\mathbb{E}(z_{k+1}|\mathcal{F}_k)\leq (1+\varsigma_k)z_k-\vartheta_k+\eta_k,~~~\forall k\in\mathbb{N}       \nonumber
\end{align}
then, $\{\vartheta_k\}_{k\in\mathbb{N}}\in\ell_+^1(\mathfrak{F})$ and $z_k$ converges to a $[0,\infty)$-valued random variable a.s.
\end{lemma}

\begin{lemma}[\cite{li2019distributed3}]\label{l6}
Let $T:\mathcal{H}\to\mathcal{H}$ be a nonexpansive operator with $Fix(T)\neq\emptyset$. Then, there holds $2\langle y-z,y-T(y)\rangle\geq \|T(y)-y\|^2$ for all $y\in\mathcal{H}$ and $z\in Fix(T)$.
\end{lemma}

\begin{lemma}\label{l7}
Under Assumptions \ref{a1}-\ref{a3}, there holds that
\begin{align}
\|x_{i,k}-\bar{x}_k\|=O(\alpha_{\lfloor\frac{k}{2}\rfloor}),~~a.s.~~\forall k\in\mathbb{N}.        \nonumber
\end{align}
\end{lemma}
{\bf Proof.}
The conclusion can be readily obtained using the same argument as in Lemma \ref{l4}, once noting that all error terms $\|b_{l,k}\alpha_k(F_{il}(\hat{x}_{i,k})-\hat{x}_{il,k})\|=O(\alpha_k)$ a.s.
\hfill\rule{2mm}{2mm}

\begin{lemma}\label{l8}
Under Assumptions \ref{a1}-\ref{a3}, one has that
\begin{align}
\Big\|\frac{\sum_{i=1}^N F_i(\hat{x}_{i,k})}{N}-\bar{x}_k\Big\|&=O(1),             \nonumber\\
\mathbb{E}\Big(|||\frac{\sum_{i=1}^N T_{i,k}}{N}-\bar{x}_k|||^2|\chi_k\Big)&=O(1),                         \nonumber\\
\mathbb{E}(|||\bar{x}_{k+1}-\bar{x}_k|||^2|\chi_k)&=O(\alpha_k^2),~~~~~a.s.         \nonumber
\end{align}
where $T_{i,k}=(T_{i1,k},\ldots,T_{im,k})$ is defined in (\ref{03}).
\end{lemma}
{\bf Proof.}
It is easy to obtain that
\begin{align}
&\Big\|\frac{\sum_{i=1}^N F_i(\hat{x}_{i,k})}{N}-\bar{x}_k\Big\|    \nonumber\\
&\leq \|\frac{\sum_{i=1}^N F_i(\hat{x}_{i,k})}{N}-F(\bar{x}_k)\|+\|F(\bar{x}_k)-\bar{x}_k\|          \nonumber\\
&\leq \frac{1}{N}\sum_{i=1}^N \|\hat{x}_{i,k}-\bar{x}_k\|+\|F(\bar{x}_k)-\bar{x}_k\|,             \label{pf16}
\end{align}
where the last inequality has exploited the nonexpansive property of $F_i$'s. In view of (\ref{pf5}), Lemma \ref{l7}, and Assumption \ref{a1}, from (\ref{pf16}) one can obtain the first conclusion of this lemma.

For the second claim of this lemma, one has that
\begin{align}
&\mathbb{E}\Big(|||\frac{\sum_{i=1}^N T_{i,k}}{N}-\bar{x}_k|||^2|\chi_k\Big)         \nonumber\\
&=\sum_{l=1}^m \frac{1}{p_l}\mathbb{E}\Big(\|\frac{\sum_{i=1}^N T_{il,k}}{N}-\bar{x}_{kl}\|^2|\chi_k\Big)    \nonumber\\
&=\sum_{l=1}^m \|\frac{\sum_{i=1}^N F_{il}(\hat{x}_{i,k})}{N}-\bar{x}_{kl}\|^2          \nonumber\\
&=\|\frac{\sum_{i=1}^N F_{i}(\hat{x}_{i,k})}{N}-\bar{x}_{k}\|^2                      \nonumber\\
&=O(1),                    \label{pf18}
\end{align}
where $\bar{x}_{kl}$ is the $l$-th component of $\bar{x}_k$, i.e., $\bar{x}_k=(\bar{x}_{k1},\ldots,\bar{x}_{km})$, and the last equality has used the first conclusion of this lemma.

Regarding the third claim, using (\ref{05a}) implies that
\begin{align}
\mathbb{E}(|||\bar{x}_{k+1}-\bar{x}_k|||^2|\chi_k)&=\alpha_k^2\mathbb{E}\Big(|||\frac{\sum_{i=1}^N T_{i,k}}{N}-\bar{x}_k|||^2|\chi_k\Big)       \nonumber\\
&=O(\alpha_k^2),                        \label{pf17}
\end{align}
where the last equality has made use of the second conclusion of this lemma.
\hfill\rule{2mm}{2mm}

Armed with the above lemma, it is now ready to give the proof of Theorem \ref{t2}.

{\em Proof of Theorem \ref{t2}:}
Note that $\bar{x}_k=(\bar{x}_{k1},\ldots,\bar{x}_{km})$ and $x^*=(x_1^*,\ldots,x_m^*)$ throughout this proof. For arbitrary $x^*\in Fix(F)$, invoking (\ref{05a}) implies that
\begin{align}
&\mathbb{E}(|||\bar{x}_{k+1}-x^*|||^2|\chi_k)            \nonumber\\
&=\mathbb{E}(|||\bar{x}_k-x^*+\alpha_k(\frac{\sum_{i=1}^N T_{i,k}}{N}-\bar{x}_k)|||^2|\chi_k)          \nonumber\\
&=|||\bar{x}_k-x^*|||^2+\alpha_k^2\mathbb{E}(|||\frac{\sum_{i=1}^N T_{i,k}}{N}-\bar{x}_k|||^2|\chi_k)           \nonumber\\
&\hspace{0.4cm}+\sum_{l=1}^m\frac{2\alpha_k}{p_l}\Big\langle \frac{\mathbb{E}(\sum_{i=1}^N T_{il,k}|\chi_k)}{N}-\bar{x}_{kl},\bar{x}_{kl}-x_l^*\Big\rangle.     \label{pf19}
\end{align}

Consider now the last two terms on the right-hand side of (\ref{pf19}). First, invoking Lemma \ref{l8} yields that for some $c_4>0$
\begin{align}
\alpha_k^2\mathbb{E}(|||\frac{\sum_{i=1}^N T_{i,k}}{N}-\bar{x}_k|||^2|\chi_k)\leq c_4\alpha_k^2\leq c_4\alpha_k\alpha_{\lfloor\frac{k}{2}\rfloor}.           \label{pf20}
\end{align}

Second, it can be obtained that
\begin{align}
&\sum_{l=1}^m\frac{2\alpha_k}{p_l}\Big\langle \frac{\mathbb{E}(\sum_{i=1}^N T_{il,k}|\chi_k)}{N}-\bar{x}_{kl},\bar{x}_{kl}-x_l^*\Big\rangle        \nonumber\\
&=2\alpha_k\sum_{l=1}^m \Big\langle \frac{\sum_{i=1}^N F_{il}(\hat{x}_{i,k})}{N}-\bar{x}_{kl},\bar{x}_{kl}-x_l^*\Big\rangle                \nonumber\\
&=\frac{2\alpha_k}{N}\sum_{i=1}^N \langle F_{i}(\hat{x}_{i,k})-F_i(\bar{x}_k),\bar{x}_k-x^*\rangle           \nonumber\\
&\hspace{0.4cm}+2\alpha_k\langle F(\bar{x}_k)-\bar{x}_{k},\bar{x}_{k}-x^*\rangle      \nonumber\\
&\leq \frac{2\alpha_k}{N}\sum_{i=1}^N \|\hat{x}_{i,k}-\bar{x}_k\|\cdot\|\bar{x}_k-x^*\|-\alpha_k\|F(\bar{x}_k)-\bar{x}_k\|^2,         \nonumber
\end{align}
where the last inequality has appealed to the Cauchy-Schwarz inequality, the nonexpansiveness of $F_i$'s, and Lemma \ref{l6}.

Furthermore, applying Young's inequality to the above inequality gives rise to
\begin{align}
&\sum_{l=1}^m\frac{2\alpha_k}{p_l}\Big\langle \frac{\mathbb{E}(\sum_{i=1}^N T_{il,k}|\chi_k)}{N}-\bar{x}_{kl},\bar{x}_{kl}-x_l^*\Big\rangle        \nonumber\\
&\leq \frac{\alpha_k}{N\alpha_{\lfloor\frac{k}{2}\rfloor}}\big(\sum_{i=1}^N \|\hat{x}_{i,k}-\bar{x}_k\|\big)^2+\frac{\alpha_k\alpha_{\lfloor\frac{k}{2}\rfloor}}{N}\|\bar{x}_k-x^*\|^2 \nonumber\\
&\hspace{0.4cm}-\alpha_k\|F(\bar{x}_k)-\bar{x}_k\|^2        \nonumber\\
&\leq\alpha_k\alpha_{\lfloor\frac{k}{2}\rfloor}\big(c_5+\frac{1}{N}\|\bar{x}_k-x^*\|^2\big)-\alpha_k\|F(\bar{x}_k)-\bar{x}_k\|^2       \label{pf21}
\end{align}
for some constant $c_5>0$, where the last inequality has employed (\ref{pf5}) and Lemma \ref{l7}.

Substituting (\ref{pf20}) and (\ref{pf21}) into (\ref{pf19}) implies that
\begin{align}
\mathbb{E}(|||\bar{x}_{k+1}-x^*|||^2|\chi_k)&\leq \Big(1+\frac{\alpha_k\alpha_{\lfloor\frac{k}{2}\rfloor}}{N}\Big)|||\bar{x}_k-x^*|||^2           \nonumber\\
&\hspace{-1.8cm}+(c_4+c_5)\alpha_k\alpha_{\lfloor\frac{k}{2}\rfloor}-\alpha_k\|F(\bar{x}_k)-\bar{x}_k\|^2,              \label{pf22}
\end{align}
which, together with Lemma \ref{l5}, results in $|||\bar{x}_k-x^*|||$ is convergent, $\bar{x}_k$ is bounded, and
\begin{align}
\sum_{k=0}^\infty\alpha_k\|F(\bar{x}_k)-\bar{x}_k\|^2<\infty,~~~a.s.              \label{pf23}
\end{align}
In light of Assumption \ref{a3}, one can obtain from (\ref{pf23}) that
\begin{align}
\mathop{\lim\inf}_{k\to\infty}\|F(\bar{x}_k)-\bar{x}_k\|^2&=0,~~~a.s.           \nonumber\\
i.e.,~~~\mathop{\lim\inf}_{k\to\infty}|||F(\bar{x}_k)-\bar{x}_k|||&=0,~~~a.s.           \label{pf24}
\end{align}

On the other hand, by (\ref{05a}), it can be concluded that
\begin{align}
&\mathbb{E}(|||F(\bar{x}_{k+1})-\bar{x}_{k+1}|||^2|\chi_k)         \nonumber\\
&=\mathbb{E}(|||F(\bar{x}_{k+1})-F(\bar{x}_k)+F(\bar{x}_k)-\bar{x}_k         \nonumber\\
&\hspace{0.4cm}-\alpha_k(\frac{\sum_{i=1}^N T_{i,k}}{N}-\bar{x}_k)|||^2|\chi_k)      \nonumber\\
&=\mathbb{E}(|||F(\bar{x}_{k+1})-F(\bar{x}_k)|||^2|\chi_k)+|||F(\bar{x}_k)-\bar{x}_k|||^2       \nonumber\\
&\hspace{0.4cm}+\alpha_k^2\mathbb{E}(|||\frac{\sum_{i=1}^N T_{i,k}}{N}-\bar{x}_k|||^2|\chi_k)      \nonumber\\
&\hspace{0.4cm}+\underbrace{\sum_{l=1}^m\frac{2}{p_l}\mathbb{E}(\langle F_{gl}(\bar{x}_{k+1})-F_{gl}(\bar{x}_k),F_{gl}(\bar{x}_k)-\bar{x}_{kl}\rangle|\chi_k)}_{=:S_6}     \nonumber\\
&\hspace{0.1cm}\underbrace{-\sum_{l=1}^m\frac{2\alpha_k}{p_l}\mathbb{E}(\langle F_{gl}(\bar{x}_{k+1})-F_{gl}(\bar{x}_k),\frac{\sum_{i=1}^N T_{il,k}}{N}-\bar{x}_{kl}\rangle|\chi_k)}_{=:S_7}  \nonumber\\
&\hspace{0.4cm}\underbrace{-\sum_{l=1}^m\frac{2\alpha_k}{p_l}\mathbb{E}(\langle F_{gl}(\bar{x}_{k})-\bar{x}_{kl},\frac{\sum_{i=1}^N T_{il,k}}{N}-\bar{x}_{kl}\rangle|\chi_k)}_{=:S_8},    \label{pf25}
\end{align}
where $F_{gl}$ is the $l$-th component of the global operator $F$, i.e., $F(x)=(F_{g1}(x),\ldots,F_{gm}(x))$ for any $x\in\mathcal{H}$.

In the sequel, consider each term in (\ref{pf25}). One has that
\begin{align}
&\mathbb{E}(|||F(\bar{x}_{k+1})-F(\bar{x}_k)|||^2|\chi_k)         \nonumber\\
&\leq \frac{1}{p_0}\mathbb{E}(\|F(\bar{x}_{k+1})-F(\bar{x}_k)\|^2|\chi_k)      \nonumber\\
&\leq \frac{1}{p_0}\mathbb{E}(\|\bar{x}_{k+1}-\bar{x}_k\|^2|\chi_k)             \nonumber\\
&\leq c_6\alpha_k^2              \label{pf26}
\end{align}
for some constant $c_6>0$, where the relationship after (\ref{09}) has been used to obtain the first inequality, the nonexpansiveness of $F$ to deduce the second inequality, and Lemma \ref{l8} to derive the last inequality.

Meanwhile, by Lemma \ref{l8}, it can be concluded that there exists a constant $c_7>0$ such that
\begin{align}
\alpha_k^2\mathbb{E}(|||\frac{\sum_{i=1}^N T_{i,k}}{N}-\bar{x}_k|||^2|\chi_k)\leq c_7\alpha_k^2.          \label{pf27}
\end{align}

For $S_6$, it can be implied that
\begin{align}
&S_6\leq\sum_{l=1}^m\frac{2}{p_l}\mathbb{E}(\|F_{gl}(\bar{x}_{k+1})-F_{gl}(\bar{x}_k)\|\|F_{gl}(\bar{x}_k)-\bar{x}_{kl}\||\chi_k)        \nonumber\\
&\leq \sum_{l=1}^m\frac{1}{p_l\alpha_k}\mathbb{E}(\|F_{gl}(\bar{x}_{k+1})-F_{gl}(\bar{x}_k)\|^2|\chi_k)           \nonumber\\
&\hspace{0.4cm}+\sum_{l=1}^m\frac{\alpha_k}{p_l}\|F_{gl}(\bar{x}_k)-\bar{x}_{kl}\|^2                        \nonumber\\
&\leq \frac{1}{p_0\alpha_k}\mathbb{E}(\|F(\bar{x}_{k+1})-F(\bar{x}_k)\|^2|\chi_k)+\frac{\alpha_k}{p_0}\|F(\bar{x}_k)-\bar{x}_k\|^2        \nonumber\\
&\leq \frac{1}{p_0\alpha_k}\mathbb{E}(\|\bar{x}_{k+1}-\bar{x}_k\|^2|\chi_k)+\frac{\alpha_k}{p_0}\|F(\bar{x}_k)-\bar{x}_k\|^2,           \nonumber
\end{align}
where the Cauchy-Schwarz inequality and nonexpansiveness of $F$ have been invoked to obtain the second and last inequalities. Applying (\ref{05a}) to the above inequality implies that
\begin{align}
&S_6\leq \frac{\alpha_k}{p_0}\mathbb{E}(\|\frac{\sum_{i=1}^N T_{i,k}}{N}-\bar{x}_k\|^2|\chi_k)+\frac{\alpha_k}{p_0}\|F(\bar{x}_k)-\bar{x}_k\|^2           \nonumber\\
&=\frac{\alpha_k}{p_0}\sum_{l=1}^m\mathbb{E}(\|\frac{\sum_{i=1}^N T_{il,k}}{N}-\bar{x}_{kl}\|^2|\chi_k)+\frac{\alpha_k}{p_0}\|F(\bar{x}_k)-\bar{x}_k\|^2    \nonumber\\
&=\frac{\alpha_k}{p_0}\sum_{l=1}^m p_l\|\frac{\sum_{i=1}^N F_{il}(\hat{x}_{i,k})}{N}-\bar{x}_{kl}\|^2+\frac{\alpha_k}{p_0}\|F(\bar{x}_k)-\bar{x}_k\|^2      \nonumber\\
&\leq \frac{\alpha_k}{p_0}\|\frac{\sum_{i=1}^N F_{i}(\hat{x}_{i,k})}{N}-\bar{x}_{k}\|^2+\frac{\alpha_k}{p_0}\|F(\bar{x}_k)-\bar{x}_k\|^2             \nonumber\\
&\leq \frac{2\alpha_k}{p_0}\|\frac{\sum_{i=1}^N F_{i}(\hat{x}_{i,k})}{N}-F(\bar{x}_k)\|^2+\frac{3\alpha_k}{p_0}\|F(\bar{x}_k)-\bar{x}_k\|^2             \nonumber\\
&\leq \frac{2\alpha_k}{Np_0}\sum_{i=1}^N\|\hat{x}_{i,k}-\bar{x}_k\|^2+\frac{3\alpha_k}{p_0}\|F(\bar{x}_k)-\bar{x}_k\|^2             \nonumber\\
&\leq c_8\alpha_k\alpha_{\lfloor\frac{k}{2}\rfloor}+\frac{3\alpha_k}{p_0}\|F(\bar{x}_k)-\bar{x}_k\|^2,                            \label{pf28}
\end{align}
where $c_8>0$ is some constant, the nonexpansiveness of $F_i$'s has been exploited to obtain the fourth inequality, and (\ref{pf5}) and Lemma \ref{l7} to have the last inequality.

For $S_7$, one has that
\begin{align}
S_7&\leq \sum_{l=1}^m\frac{2\alpha_k}{p_l}\mathbb{E}(\|F_{gl}(\bar{x}_{k+1})-F_{gl}(\bar{x}_k)\|       \nonumber\\
&\hspace{0.4cm}\cdot\|\frac{\sum_{i=1}^N T_{il,k}}{N}-\bar{x}_{kl}\||\chi_k),      \nonumber
\end{align}
where, invoking Assumption \ref{a1}, it is easy to verify that
\begin{align}
\|\frac{\sum_{i=1}^N T_{il,k}}{N}-\bar{x}_{kl}\|&=\|\frac{\sum_{i=1}^N b_{il.k}(F_{il}(\hat{x}_{i,k})-\hat{x}_{il,k})}{N}\|      \nonumber\\
&\leq\frac{1}{N}\sum_{i=1}^N \|F_{il}(\hat{x}_{i,k})-\hat{x}_{il,k}\|        \nonumber\\
&\leq B,~~~~~~~~~~~a.s.         \nonumber
\end{align}
which further implies that
\begin{align}
&S_7\leq \sum_{l=1}^m\frac{2B\alpha_k}{p_l}\mathbb{E}(\|F_{gl}(\bar{x}_{k+1})-F_{gl}(\bar{x}_k)\||\chi_k)      \nonumber\\
&\leq \frac{2\sqrt{m}B\alpha_k}{p_0}\mathbb{E}(\|F(\bar{x}_{k+1})-F(\bar{x}_k)\||\chi_k)      \nonumber\\
&\leq \frac{2\sqrt{m}B\alpha_k}{p_0}\mathbb{E}(\|\bar{x}_{k+1}-\bar{x}_k\||\chi_k)      \nonumber\\
&\leq \frac{2\sqrt{m}B\alpha_k}{p_0}\sqrt{\mathbb{E}(\|\bar{x}_{k+1}-\bar{x}_k\|^2|\chi_k)}      \nonumber\\
&\leq c_9\alpha_k^2                  \label{pf29}
\end{align}
for some constant $c_9>0$, where the fact that $\sum_{l=1}^m\|z_l\|\leq \sqrt{m}\|z\|$ for any $z=(z_1,\ldots,z_m)$ has been used to get the second inequality, the nonexpansiveness of $F$ to obtain the third inequality, Jensen's inequality to deduce the fourth inequality, and Lemma \ref{l8} to derive the last inequality.

For $S_8$, it can be concluded that
\begin{align}
S_8&=-\sum_{l=1}^m\frac{2\alpha_k}{p_l}\langle F_{gl}(\bar{x}_k)-\bar{x}_{kl},\mathbb{E}(\frac{\sum_{i=1}^N T_{il,k}}{N}-\bar{x}_{kl}|\chi_k)\rangle      \nonumber\\
&=-\sum_{l=1}^m 2\alpha_k \langle F_{gl}(\bar{x}_k)-\bar{x}_{kl},\frac{\sum_{i=1}^N F_{il}(\hat{x}_{i,k})}{N}-\bar{x}_{kl}\rangle      \nonumber\\
&=-2\alpha_k \|F(\bar{x}_k)-\bar{x}_{k}\|^2          \nonumber\\
&\hspace{0.4cm}+2\alpha_k\sum_{l=1}^m \langle F_{gl}(\bar{x}_k)-\bar{x}_{kl},F_{gl}(\bar{x}_k)-\frac{\sum_{i=1}^N F_{il}(\hat{x}_{i,k})}{N}\rangle      \nonumber\\
&\leq 2\alpha_k\sum_{l=1}^m \|F_{gl}(\bar{x}_k)-\bar{x}_{kl}\|\|F_{gl}(\bar{x}_k)-\frac{\sum_{i=1}^N F_{il}(\hat{x}_{i,k})}{N}\|          \nonumber\\
&\leq 2B\alpha_k\sum_{l=1}^m \|F_{gl}(\bar{x}_k)-\frac{\sum_{i=1}^N F_{il}(\hat{x}_{i,k})}{N}\|          \nonumber\\
&\leq 2B\sqrt{m}\alpha_k \|F(\bar{x}_k)-\frac{\sum_{i=1}^N F_{i}(\hat{x}_{i,k})}{N}\|              \nonumber\\
&\leq \frac{2B\sqrt{m}\alpha_k}{N} \sum_{i=1}^N\|F_i(\bar{x}_k)-F_{i}(\hat{x}_{i,k})\|          \nonumber\\
&\leq \frac{2B\sqrt{m}\alpha_k}{N} \sum_{i=1}^N\|\bar{x}_k-\hat{x}_{i,k}\|          \nonumber\\
&\leq c_{10}\alpha_k\alpha_{\lfloor\frac{k}{2}\rfloor}                            \label{pf30}
\end{align}
for some constant $c_{10}>0$, where Assumption \ref{a1} has been appealed to obtain the second inequality, the nonexpansiveness of $F_i$'s to deduce the fifth inequality, and (\ref{pf5}) and Lemma \ref{l7} to derive the last inequality.

At this position, substituting (\ref{pf26})-(\ref{pf30}) into (\ref{pf25}) yields that
\begin{align}
\mathbb{E}(|||F(\bar{x}_{k+1})-\bar{x}_{k+1}|||^2|\chi_k)&\leq |||F(\bar{x}_{k})-\bar{x}_{k}|||^2           \nonumber\\
&\hspace{-2.1cm}+c_1'\alpha_k\alpha_{\lfloor\frac{k}{2}\rfloor}+\frac{3\alpha_k}{p_0}\|F(\bar{x}_k)-\bar{x}_k\|^2,        \label{pf31}
\end{align}
where the fact $\alpha_k^2\leq \alpha_k\alpha_{\lfloor\frac{k}{2}\rfloor}$ has been utilized, and $c_1':=\sum_{d=6}^{10}c_d$.

Now, applying Lemma \ref{l5} to (\ref{pf31}), in conjunction with (\ref{pf23}) and Assumption \ref{a3}, one can conclude that $|||F(\bar{x}_k)-\bar{x}_k|||^2$ is convergent a.s., which, together with (\ref{pf24}), further implies that
\begin{align}
\lim_{k\to\infty}|||F(\bar{x}_k)-\bar{x}_k|||=0,~~~a.s.      \label{pf32}
\end{align}

Finally, following the same argument as that after (\ref{pf14}), Theorem \ref{t2} can be claimed. This completes the proof.
\hfill\rule{2mm}{2mm}

\subsection{Proof of Proposition \ref{p1}}\label{s7.3}

For brevity, $P_{X_i}$ is still denoted as $F_i$ for $i\in[N]$ throughout this proof.

First, the D-KM iteration (\ref{x6}) can be written as
\begin{align}
x_{i,k+1}=(1-\alpha_k)\hat{x}_{i,k}+\alpha_k F_i(\hat{x}_{i,k}).          \label{pf33}
\end{align}

Note that $|1-\alpha_k|\leq 1$ and $F_i$'s are bounded due to the compactness of $X_i$. Hence, similar to Lemma \ref{l4}, one can obtain that (\ref{pf2}) holds. Now, invoking the same argument as in (\ref{pf4})-(\ref{pf7}) yields that $\bar{x}_k$ is bounded, which, together with (\ref{pf2}), further implies that $x_{i,k}$'s are all bounded. Notice that $F_i$'s are bounded. As a result, Assumption \ref{a1} holds for (\ref{pf6}).

As for the D-BKM iteration (\ref{x7}), which amounts to
\begin{align}
x_{i,k+1}=(I-\alpha_k\Gamma_k)\hat{x}_{i,k}+\alpha_k\Gamma_k F_i(\hat{x}_{i,k}),            \label{pf34}
\end{align}
where $\Gamma_k:=diag\{b_{1,k}I,\ldots,b_{m,k}I\}$.

Keeping in mind (\ref{pf34}) and $\|1-\alpha_k\Gamma_k\|\leq 1$, appealing to the same reasoning as the above for (\ref{pf33}) gives rise to that (\ref{pf2}) holds a.s.

Subsequently, for $x^*\in Fix(F)$, in view of (\ref{05a}), it can be deduced that
\begin{align}
&\mathbb{E}(\|\bar{x}_{k+1}-x^*\|^2|\chi_k)         \nonumber\\
&=\mathbb{E}(\|(1-\alpha_k)(\bar{x}_k-x^*)+\alpha_k(\frac{\sum_{i=1}^N T_{i,k}}{N}-x^*)\|^2|\chi_k)     \nonumber\\
&\leq (1-\alpha_k)\|\bar{x}_k-x^*\|^2+\alpha_k\sum_{l=1}^m\mathbb{E}(\|\frac{\sum_{i=1}^N T_{il,k}}{N}-x_l^*\|^2|\chi_k)        \nonumber\\
&=(1-\alpha_k)\|\bar{x}_k-x^*\|^2+\alpha_k(1-\frac{1}{m})\sum_{l=1}^m \|\bar{x}_{kl}-x_l^*\|^2          \nonumber\\
&\hspace{0.4cm}+\frac{\alpha_k}{m}\sum_{l=1}^m\|\frac{\sum_{i=1}^N F_{il}(\hat{x}_{i,k})}{N}-x_l^*\|^2        \nonumber
\end{align}
\begin{align}
&=(1-\alpha_k)\|\bar{x}_k-x^*\|^2+\alpha_k(1-\frac{1}{m}) \|\bar{x}_{k}-x^*\|^2              \nonumber\\
&\hspace{0.4cm}+\frac{\alpha_k}{m}\|\frac{\sum_{i=1}^N F_{i}(\hat{x}_{i,k})}{N}-x^*\|^2,       \label{pf35}
\end{align}
where Lemma \ref{l2} and $p_l=1/m,i\in[m]$ have been applied to obtain the first inequality and second equality, respectively.

Regarding the last term in (\ref{pf35}), one has that
\begin{align}
&\|\frac{\sum_{i=1}^N F_{i}(\hat{x}_{i,k})}{N}-x^*\|^2         \nonumber\\
&=\|\frac{\sum_{i=1}^N F_{i}(\hat{x}_{i,k})}{N}-F(\bar{x}_k)+F(\bar{x}_k)-x^*\|^2         \nonumber\\
&=\|\frac{\sum_{i=1}^N F_{i}(\hat{x}_{i,k})}{N}-F(\bar{x}_k)\|^2+\|F(\bar{x}_k)-x^*\|^2           \nonumber\\
&\hspace{0.4cm}+2\langle \frac{\sum_{i=1}^N F_{i}(\hat{x}_{i,k})}{N}-F(\bar{x}_k),F(\bar{x}_k)-x^*\|^2\rangle    \nonumber\\
&\leq \frac{1}{N^2}(\sum_{i=1}^N\|\hat{x}_{i,k}-\bar{x}_k\|)^2+\|\bar{x}_k-x^*\|^2              \nonumber\\
&\hspace{0.4cm}+\frac{2}{N}\sum_{i=1}^N \|\hat{x}_{i,k}-\bar{x}_k\|\cdot\|\bar{x}_k-x^*\|          \nonumber\\
&\leq \frac{1}{N^2}(\sum_{i=1}^N\|\hat{x}_{i,k}-\bar{x}_k\|)^2+\|\bar{x}_k-x^*\|^2              \nonumber\\
&\hspace{0.4cm}+\frac{\alpha_{\lfloor\frac{k}{2}\rfloor}}{N}\|\bar{x}_k-x^*\|^2+\frac{(\sum_{i=1}^N\|\hat{x}_{i,k}-\bar{x}_k\|)^2}{N\alpha_{\lfloor\frac{k}{2}\rfloor}},        \label{pf36}
\end{align}
where the Cauchy-Schwarz inequality and the nonexpansiveness of $F$ and $F_i$'s have been invoked to obtained the first inequality, and Young's inequality has been utilized to get the last inequality.

In light of (\ref{pf5}) and (\ref{pf2}), plugging (\ref{pf36}) into (\ref{pf35}) yields that
\begin{align}
&\mathbb{E}(\|\bar{x}_{k+1}-x^*\|^2|\chi_k)         \nonumber\\
&\leq (1+\frac{\alpha_k\alpha_{\lfloor\frac{k}{2}\rfloor}}{mN})\|\bar{x}_k-x^*\|^2+\zeta_1\alpha_k\alpha_{\lfloor\frac{k}{2}\rfloor},       \label{pf37}
\end{align}
where $\zeta_1>0$ is some constant.

Invoking Lemma \ref{l5} and Assumption \ref{a3}, one can claim from (\ref{pf37}) that $\bar{x}_k$ is bounded a.s., which in combination with (\ref{pf2}) ensures the boundedness of $x_{i,k}$'s a.s. The conclusion can be then drawn.
\hfill\rule{2mm}{2mm}

\subsection{Proof of Proposition \ref{p2}}\label{s7.4}

In this case, the D-KM iteration can be written as
\begin{align}
x_{i,k+1}=\Omega_{i,k}\hat{x}_{i,k}+\alpha_k\theta r_i,              \label{pf38}
\end{align}
where $\Omega_{i,k}:=I-\alpha_k\theta R_i$. It should be noted that $\|\Omega_{i,k}\|\leq 1$ for all $i\in[N],k\in\mathbb{N}$ due to $\theta\in(0,2/\lambda_M]$ and $\alpha_k\in(0,1]$.

Moreover, (\ref{pf38}) can be written in a compact form
\begin{align}
x_{k+1}=T_k(x_k),              \label{pf39}
\end{align}
where $x_k=col(x_{1,k},\ldots,x_{N,k})$ and $T_k$ is defined as
\begin{align}
T_k:x\mapsto \Omega_k (A_k\otimes I_n) x+\alpha_k\theta \tilde{r},
\end{align}
with $\Omega_k:=diag\{I-\alpha_k\theta R_1,\ldots,I-\alpha_k\theta R_N\}$ and $\tilde{r}:=col(r_1,\ldots,r_N)$.

It is easy to verify that for any $x,y\in\mathbb{R}^{nN}$
\begin{align}
\|T_k x-T_ky\|&=\|\Omega_k(A_k\otimes I_n)x-\Omega_k(A_k\otimes I_n)y\|       \nonumber\\
&\leq \|\Omega_k\|\|(A_k\otimes I_n)x-(A_k\otimes I_n)y\|                     \nonumber\\
&\leq \|(A_k\otimes I_n)x-(A_k\otimes I_n)y\|                           \nonumber\\
&\leq \|A_k\otimes I_n\|\|x-y\|                             \nonumber\\
&= \|x-y\|,                              \label{pf40}
\end{align}
where $\|\Omega_k\|\leq 1$ and $\|A_k\otimes I_n\|=1$ have been employed. Therefore, $T_k$ is nonexpansive.

Now, for any $x_k^*\in Fix(T_k)$, invoking (\ref{pf39}) yields that
\begin{align}
\|x_{k+1}-x_k^*\|&=\|T_kx_k-T_kx_k^*\|\leq \|x_k-x_k^*\|,                              \label{pf41}
\end{align}
implying that $x_k$ and thus $x_{i,k}$'s are bounded. Consequently, Assumption \ref{a1} holds. This ends the proof.
\hfill\rule{2mm}{2mm}


\begin{thebibliography}{10}

\bibitem{alaviani2018distributed}
Alaviani, S.~S., \& Elia, N. (2018).
\newblock A distributed algorithm for solving linear algebraic equations over
  random networks.
\newblock In {\em Proceedings of 57th Conference on Decision and Control},
  pp. 83--88, Miami Beach, FL, USA.

\bibitem{alaviani2019distributed}
Alaviani, S.~S., \& Elia, N. (2019).
\newblock Distributed multi-agent convex optimization over random digraphs.
\newblock {\em IEEE Transactions on Automatic Control}, in press, doi:
  10.1109/TAC.2019.2937499.

\bibitem{bauschke2017convex}
Bauschke, H.~H., \& Combettes, P.~L. (2017).
\newblock {\em {Convex analysis and monotone operator theory in Hilbert
  spaces}, 2nd ed}.
\newblock Springer, New York.

\bibitem{borwein2017convergence}
Borwein, J.~M., Li, G., \& Tam, M.~K. (2017).
\newblock Convergence rate analysis for averaged fixed point iterations in
  common fixed point problems.
\newblock {\em SIAM Journal on Optimization}, {\em 27}(1), 1--33.

\bibitem{bravo2018rates}
Bravo, M., Cominetti, R., \& Pavez-Sign{\'e}, M. (2019).
\newblock {Rates of convergence for inexact Krasnosel'ski\u{\i}-Mann iterations
  in Banach spaces}.
\newblock {\em Mathematical Programming}, {\em 175}(1-2), 241--262.

\bibitem{cegielski2012iterative}
Cegielski, A. (2012).
\newblock {\em Iterative methods for fixed point problems in Hilbert spaces},
  Vol. 2057.
\newblock Springer, Heidelberg.

\bibitem{cegielski2015application}
Cegielski, A. (2015).
\newblock Application of quasi-nonexpansive operators to an iterative method
  for variational inequality.
\newblock {\em SIAM Journal on Optimization}, {\em 25}(4), 2165--2181.

\bibitem{cominetti2014rate}
Cominetti, R., Soto, J.~A., \& Vaisman, J. (2014).
\newblock {On the rate of convergence of Krasnosel'ski\u{\i}-Mann iterations
  and their connection with sums of Bernoullis}.
\newblock {\em Israel Journal of Mathematics}, {\em 199}(2), 757--772.

\bibitem{dall2019convergence}
Dall'Anese, E., Simonetto, A., \& Bernstein, A. (2019).
\newblock On the convergence of the inexact running Krasnosel'ski\u{\i}-Mann
  method.
\newblock {\em IEEE Control Systems Letters}, {\em 3}(3), 613--618.

\bibitem{fullmer2016asynchronous}
Fullmer, D., Liu, J., \& Morse, A.~S. (2016).
\newblock An asynchronous distributed algorithm for computing a common fixed
  point of a family of paracontractions.
\newblock In {\em Proceedings of 55th Conference on Decision and Control},
  pp. 2620--2625, Las Vegas, USA.

\bibitem{fullmer2018distributed}
Fullmer, D., \& Morse, A.~S. (2018).
\newblock A distributed algorithm for computing a common fixed point of a
  finite family of paracontractions.
\newblock {\em IEEE Transactions on Automatic Control}, {\em 63}(9), 2833--2843.

\bibitem{iiduka2016convergence}
Iiduka, H. (2016).
\newblock Convergence analysis of iterative methods for nonsmooth convex
  optimization over fixed point sets of quasi-nonexpansive mappings.
\newblock {\em Mathematical Programming}, {\em 159}(1-2), 509--538.

\bibitem{kanzow2017generalized}
Kanzow, C., \& Shehu, Y. (2017).
\newblock {Generalized Krasnosel'ski\u{\i}-Mann-type iterations for
  nonexpansive mappings in Hilbert spaces}.
\newblock {\em Computational Optimization and Applications}, {\em 67}(3), 595--620.

\bibitem{krasnosel1955two}
Krasnosel'ski\u{\i}, M.~A. (1955).
\newblock Two comments on the method of successive approximations.
\newblock {\em Uspekhi Matematicheskikh Nauk}, {\em 10}, 123--127.

\bibitem{kruger2018set}
Kruger, A. Y., Luke, D.~R., \& Thao, N.~H. (2018).
\newblock Set regularities and feasibility problems.
\newblock {\em Mathematical Programming}, {\em 168}(1-2), 279--311.

\bibitem{li2016consensus}
Li, X., Chen, M.~Z.~Q., Su, H., \& Li, C. (2016).
\newblock Consensus networks with switching topology and time-delays over
  finite fields.
\newblock {\em Automatica}, {\em 68}, 39--43.

\bibitem{li2019distributed3}
Li, X., \& Feng, G. (2019).
\newblock Distributed algorithms for computing a common fixed point of a group
  of nonexpansive operators.
\newblock arXiv preprint arXiv:1902.02481.

\bibitem{li2019distributed}
Li, X., Xie, L., \& Hong, Y. (2019).
\newblock Distributed continuous-time nonsmooth convex optimization with
  coupled inequality constraints.
\newblock {\em IEEE Transactions on Control of Network Systems}, in press, doi:
  10.1109/TCNS.2019.2915626.

\bibitem{liang2016convergence}
Liang, J., Fadili, J., \& Peyr{\'e}, G. (2016).
\newblock Convergence rates with inexact non-expansive operators.
\newblock {\em Mathematical Programming}, {\em 159}(1-2), 403--434.

\bibitem{liang2019distributed}
Liang, S., Wang, L., \& Yin, G. (2019).
\newblock Distributed quasi-monotone subgradient algorithm for nonsmooth convex
  optimization over directed graphs.
\newblock {\em Automatica}, {\em 101}, 175--181.

\bibitem{lin2018multiagent}
Lin, P., Ren, W., Wang, H., \& Al-Saggaf, U.~M. (2019).
\newblock {Multi-agent rendezvous with shortest distance to convex regions with
  empty intersection: algorithms and experiments}.
\newblock {\em IEEE Transactions on Cybernetics}, {\em 49}(3), 1026--1034.

\bibitem{liu2017distributed}
Liu, J., Fullmer, D., Nedi{\'c}, A., Ba{\c{s}}ar, T., \& Morse, A.~S. (2017).
\newblock A distributed algorithm for computing a common fixed point of a
  family of strongly quasi-nonexpansive maps.
\newblock In {\em Proceedings of American Control Conference}, pp. 686--690,
  Seattle, USA.

\bibitem{liu2017convergence}
Liu, S., Qiu, Z., \& Xie, L. (2017).
\newblock Convergence rate analysis of distributed optimization with projected
  subgradient algorithm.
\newblock {\em Automatica}, {\em 83}, 162--169.

\bibitem{lou2016distributed}
Lou, Y., Hong, Y., \& Wang, S. (2016).
\newblock Distributed continuous-time approximate projection protocols for
  shortest distance optimization problems.
\newblock {\em Automatica}, {\em 69}, 289--297.

\bibitem{mann1953mean}
Mann, W.~R. (1953).
\newblock Mean value methods in iteration.
\newblock {\em Proceedings of the American Mathematical Society},
  {\em 4}(3), 506--510.

\bibitem{mansoori2019a}
Mansoori, F., \& Wei, E. (2019).
\newblock {A fast distributed asynchronous Newton-based optimization
  algorithm}.
\newblock {\em IEEE Transactions on Automatic Control}, in press, doi:
  10.1109/TAC.2019.2907711.

\bibitem{mateos2017distributed}
Mateos-N{\'u}nez, D., \& Cort{\'e}s, J. (2017).
\newblock {Distributed saddle-point subgradient algorithms with Laplacian
  averaging}.
\newblock {\em IEEE Transactions on Automatic Control}, {\em 62}(6), 2720--2735.

\bibitem{matsushita2017convergence}
Matsushita, S. (2017).
\newblock {On the convergence rate of the Krasnosel'ski\u{\i}-Mann iteration}.
\newblock {\em Bulletin of the Australian Mathematical Society},
  {\em 96}(1), 162--170.

\bibitem{meng2017adaptive}
Meng, M., Liu, L., \& Feng, G. (2017).
\newblock {Adaptive output regulation of heterogeneous multi-agent systems
  under Markovian switching topologies}.
\newblock {\em IEEE Transactions on Cybernetics}, {\em 48}(10), 2962--2971.

\bibitem{mou2015distributed}
Mou, S., Liu, J., \& Morse, A.~S. (2015).
\newblock A distributed algorithm for solving a linear algebraic equation.
\newblock {\em IEEE Transactions on Automatic Control}, {\em 60}(11), 2863--2878.

\bibitem{necoara2018randomized}
Necoara, I., Richt{\'a}rik, P., \& Patrascu, A. (2018).
\newblock {Randomized projection methods for convex feasibility problems:
  conditioning and convergence rates}.
\newblock arXiv preprint arXiv:1801.04873.

\bibitem{nedic2009distributed}
Nedi\'{c}, A., \& Ozdaglar, A. (2009).
\newblock Distributed subgradient methods for multi-agent optimization.
\newblock {\em IEEE Transactions on Automatic Control}, {\em 54}(1), 48--61.

\bibitem{nedic2018multi}
Nedi{\'c}, A., Pang, J., Scutari, G., \& Sun, Y. (2018).
\newblock {\em Multi-agent optimization: Cetraro, Italy 2014}, Vol. 2224.
\newblock Springer.

\bibitem{olfati2007consensus}
Olfati-Saber, R., Fax, J.~A., \& Murray, R.~M. (2007).
\newblock Consensus and cooperation in networked multi-agent systems.
\newblock {\em Proceedings of the IEEE}, {\em 95}(1), 215--233.

\bibitem{reich1979weak}
Reich, S. (1979).
\newblock {Weak convergence theorems for nonexpansive mappings in Banach
  spaces}.
\newblock {\em Journal of Mathematical Analysis and Applications},
  {\em 67}(2), 274--276.

\bibitem{ren2010distributed}
Ren, W., \& Cao, Y. (2010).
\newblock {\em {Distributed coordination of multi-agent networks: emergent
  problems, models, and issues}}.
\newblock London, U.K.: Springer-Verlag.

\bibitem{robbins1971martin}
Robbins, H., \& Siegmund, D. (1971).
\newblock A convergence theorem for nonnegative almost supermartingales and
  some applications.
\newblock In {\em Proceedings of Optimizing Methods in Statistics}, pp. 233--257, Ohio, USA.

\bibitem{shehu2018convergence}
Shehu, Y. (2018).
\newblock {Convergence rate analysis of inertial Krasnosel'ski\u{\i}-Mann type
  iteration with applications}.
\newblock {\em Numerical Functional Analysis and Optimization},
  {\em 39}(10), 1077--1091.

\bibitem{themelis2019supermann}
Themelis, A., \& Patrinos, P. (2019).
\newblock {SuperMann: a superlinearly convergent algorithm for finding fixed
  points of nonexpansive operators}.
\newblock {\em IEEE Transactions on Automatic Control}, in press, doi:
  10.1109/TAC.2019.2906393.

\bibitem{wang2019distributed}
Wang, P., Ren, W., \& Duan, Z. (2019).
\newblock Distributed algorithm to solve a system of linear equations with
  unique or multiple solutions from arbitrary initializations.
\newblock {\em IEEE Transactions on Control of Network Systems}, {\em 6}(1), 82--93.

\bibitem{wang2017further}
Wang, X., Mou, S., \& Sun, D. (2017).
\newblock Further discussions on a distributed algorithm for solving linear
  algebra equations.
\newblock In {\em Proceedings of American Control Conference}, pp. 4274--4278, Seattle, USA.

\bibitem{wang2019scalable}
Wang, X., Mou, S., \& Anderson, B.~D.~O. (2019).
\newblock Scalable, distributed algorithms for solving linear equations via
  double-layered networks.
\newblock {\em IEEE Transactions on Automatic Control}, in press, doi:
  10.1109/TAC.2019.2919101.

\bibitem{xie2018distributed}
Xie, P., You, K., Tempo, R., Song, S., \& Wu, C. (2018).
\newblock Distributed convex optimization with inequality constraints over
  time-varying unbalanced digraphs.
\newblock {\em IEEE Transactions on Automatic Control}, {\em 63}(12), 4331--4337.

\bibitem{yang2019survey}
Yang, T., Yi, X., Wu, J., Yuan, Y., Wu, D., Meng, Z., Hong, Y., Wang, H., Lin, Z., \& Johansson, K. H. (2019).
\newblock A survey of distributed optimization.
\newblock {\em Annual Reviews in Control}, {\em 47}, 278--305.

\bibitem{zeng2019generalized}
Zeng, X., Chen, J., Liang, S., \& Hong, Y. (2019).
\newblock {Generalized Nash equilibrium seeking strategy for distributed
  nonsmooth multi-cluster game}.
\newblock {\em Automatica}, {\em 103}, 20--26.

\bibitem{zhang2015distributed}
Zhang, Y., Lou, Y., Hong, Y., \& Xie, L. (2015).
\newblock Distributed projection-based algorithms for source localization in
  wireless sensor networks.
\newblock {\em IEEE Transactions on Wireless Communications}, {\em 14}(6), 3131--3142.

\end{thebibliography}
\end{document}